\documentclass{amsart}

\usepackage[monochrome]{xcolor}

\usepackage[utf8]{inputenc}
\usepackage[TS1,T1]{fontenc}
\usepackage{lmodern} 
\usepackage{amssymb} 
\usepackage{amsmath} 
\usepackage{bm}
\usepackage{figlatex}
\usepackage{graphics}
\graphicspath{{figures/}}
\numberwithin{equation}{section}

\newtheorem{theorem}{Theorem}[section]
\newtheorem{proposition}[theorem]{Proposition}
\newtheorem{lemma}[theorem]{Lemma}
\newtheorem{corollary}[theorem]{Corollary}

\theoremstyle{definition}
\newtheorem{definition}[theorem]{Definition}
\newtheorem{remark}[theorem]{Remark}

\usepackage{enumitem} 
\newlist{myenum}{enumerate}{1}  
\setlist[myenum,1]{label=(\roman*),%
                  ref=(\roman*)}
\newenvironment{mydescription}{\begin{enumerate}}{\end{enumerate}}

\usepackage{hyperref}
\hypersetup{
  colorlinks=true,    %
  breaklinks=true,    %
   citecolor= blue,
   linkcolor= red,    %
   bookmarksopen=true,         %
   pdftitle={On Triples...},  %
   pdfauthor={A. Parreau},          %
    }
\newcommand{\pinfty}{{+\infty}}         
\newcommand{\minfty}{{-\infty}}
\newcommand{\eps}{\varepsilon}
\newcommand{\ov}{\overline}
\newcommand{\ra}{\rightarrow} 
\newcommand{\ovra}{\overrightarrow} 
\newcommand{\isomto}{\stackrel{\sim}{\to}}
\newcommand{\abs}[1]{\left\lvert #1 \right\rvert} 
\newcommand{\RR}{\mathbb{R}}
\newcommand{\NN}{\mathbb{N}}
\newcommand{\ZZ}{\mathbb{Z}}
\newcommand{\CC}{\mathbb{C}}

\newcommand{\pospart}[1]{#1^+} 
\newcommand{\negpart}[1]{#1^-} 

\DeclareMathOperator{\GL}{GL}        
\DeclareMathOperator{\SL}{SL}        
\DeclareMathOperator{\PGL}{PGL}

\DeclareMathOperator{\diag}{diag}

\newcommand{\bord}{\partial}    
\newcommand{\bordinf}{\bord_\infty} 

\newcommand{\HH}{\mathbb{H}}

\newcommand{\ii}{i} 
\newcommand{\jj}{j} 
\newcommand{\kk}{k} 

\newcommand{\N}{N} 
\newcommand{\M}{M} 
\newcommand{\indN}{i} 
\newcommand{\indNp}{j} 
 \newcommand{\indNpp}{k} 
\newcommand{\indTriRap}{m} 

\newcommand{\n}{n} 

\newcommand{\sR}{s} 
\newcommand{\tR}{t} 

\newcommand{\KK}{\mathbb{K}}
\newcommand{\aK}{a} 
\newcommand{\bK}{b} 

\newcommand{\G}{G}        
\newcommand{\g}{g}  
\newcommand{\aG}{a} 

\newcommand{\V}{V} 
 
\newcommand{\eVb}{\mathbf{v}} 
\newcommand{\eV}{v} 
\newcommand{\normV}{\eta}
\newcommand{\PP}{\mathbb{P}} 
\newcommand{\U}{U} 
\newcommand{\independent}{independent} 
\newcommand{\Frame}{\mathcal{F}}

\newcommand{\PV}{\PP(\V)} 
\newcommand{\classPV}[1]{[#1]} 

\DeclareMathOperator{\AllFlags}{Flags}

\newcommand{\p}{p} 
\newcommand{\q}{q} 
\newcommand{\D}{D} 
\newcommand{\F}{F}

\newcommand{\pt}{{\tilde \p}} 
                              
\newcommand{\Dt}{{\tilde \D}}

\newcommand{\nondegenerateQ}{nondegenerate}
\DeclareMathOperator{\Bir}{\mathbf{b}}

\DeclareMathOperator{\geombir}{\bm{\beta}}

\newcommand{\Aa}{\mathbb{A}} 

\newcommand{\vRN}{\lambda} 

\newcommand{\classA}[1]{[{#1}]} 
\newcommand{\classAcoordsll}[2]{\classA{(#1,#2)}}

\newcommand{\eRRNo}{\eps} 

\newcommand{\vA}{\lambda} 

\newcommand{\goth}[1]{\mathfrak{#1}}
\newcommand{\Cc}{{\goth C}} 
\newcommand{\Cb}{\overline{\Cc}} 

\newcommand{\W}{W}    
\newcommand{\Sym}{{\mathfrak S}}        
\newcommand{\Waff}{W_{aff}}    

\newcommand{\stypep}{{1}} 
 
\newcommand{\stypeH}{{{\N-1}}} 

\newcommand{\rac}[1]{\varphi_{#1}} 
\newcommand{\RacSimples}{\Lambda} 
\newcommand{\deucl}{d}

\newcommand{\B}{\mathcal{B}}

\DeclareMathOperator{\CAT}{CAT}
\newcommand{\Es}{X} 

\newcommand{\Y}{Y} 

\newcommand{\f}{f} 

\newcommand{\xE}{x}  
\newcommand{\xEp}{x'}  
\newcommand{\yE}{y} 
\newcommand{\zE}{z}

\newcommand{\bordinfEs}{\bordinf\Es} 

\newcommand{\ray}{r} 

\newcommand{\xT}{x}  
\newcommand{\yT}{y}  

\newcommand{\App}{A}   
\newcommand{\Ch}{C}   

 \newcommand{\SectE}{S} 

 \newcommand{\SW}{S} 

 \DeclareMathOperator{\TangS}{\Sigma}

\newcommand{\vE}{v}

\newcommand{\ch}{c}
\newcommand{\chp}{c_+}
\newcommand{\chm}{c_-}   
\newcommand{\Star}{St}   
\newcommand{\ret}{r} 
\newcommand{\Bus}[3]{B_{#1}(#2, #3)} 
\newcommand{\Buso}{B}

\newcommand{\xim}{{\xi^-}} 
\newcommand{\xip}{{\xi^+}} 
\newcommand{\ximp}{{\xim\xip}} 
\newcommand{\Esximp}{F_\ximp}

\newcommand{\dxi}{d_\xi} 
\DeclareMathOperator{\proj}{proj}
\newcommand{\projP}{\proj} 
\newcommand{\projE}{\pi}

\newcommand{\centre}{c} 

\newcommand{\Tau}{T}

\newcommand{\Aij}{\App_{{\ii\jj}}} 
\newcommand{\Ajk}{\App_{{\jj\kk}}} 
\newcommand{\Aki}{\App_{{\kk\ii}}} 
 
\newcommand{\Aik}{\App_{\ii\kk}} 
\newcommand{\Ap}{\App_{\p}} 
\newcommand{\AD}{\App_{\D}} 
\newcommand{\yTauE}{y}  
\newcommand{\vstE}{x}  
\newcommand{\typeflat}{``flat''} 
\newcommand{\typetripod}{``tripod''} 

 \newcommand{\nondegenerateTF}{nondegenerate}
 \newcommand{\NondegenerateTF}{Nondegenerate}
 \newcommand{\ND}{ND}
 \newcommand{\genericTF}{generic}

\DeclareMathOperator{\Tri}{Tri} 
\newcommand{\Z}{Z} 

\DeclareMathOperator{\geomtri}{tri} 
\newcommand{\gZ}{z}

\begin{document}

\title{On triples of ideal chambers in  $A_2$-buildings}

\author[A. Parreau]{Anne Parreau}
\address{%
  Universit{\'e} Grenoble Alpes et CNRS\\
  Institut Fourier,\\ 
  Grenoble\\
  France
}
\email{Anne.Parreau@univ-grenoble-alpes.fr}



\begin{abstract}
  We investigate the geometry in a real Euclidean building $X$ of
  type $A_2$ of some simple configurations in the associated
  projective plane at infinity $\mathbb{P}$, seen as ideal configurations in
  $X$, and relate it with the projective invariants (from the cross
  ratio on $\mathbb{P}$).  In particular we establish a geometric
  classification of generic triples of ideal chambers of $X$ and
  relate it with the triple ratio of triples of flags.
\end{abstract}

\maketitle

\section*{Introduction}
The triples of objects in the boundaries of 
geometric spaces $\Es$ 
are basic tools, 
for example in the study of surface group representations.
For instance, in the case where $\Es=\HH^2$,
ideal triples of points may be used to define
the notion of Euler class  \cite{Gol80}, and
Penner-Thurston shear coordinates on the Teichm{\"u}ller space.
In the case where  $\Es=\HH^2_\CC$, 
the  ideal triples
are classified by   Cartan's angular invariant, 
see for example \cite[\S  7.1]{GolBook}, 
and they may be for instance used  to
define Toledo's invariant and maximal representations,  
see \cite{Tol89}. 
See for instance \cite{ClNe06,BIW10} for generalization to higher rank
Hermitian symmetric spaces $\Es$, and   triples in their
Shilov boundary.

For higher rank  symmetric spaces $\Es$ of type $A_{\N-1}$,
corresponding to the group $\PGL_\N(\RR)$, 
ideal configurations in $\Es$ 
may be seen as configurations in the
projective space $\PP=\PP(\RR^\N)$.
In particular, ideal chambers of $\Es$ 
correspond to complete flags in $\PP$,
and opposite pairs of flags (or generic $\N$-tuples of points)
in $\PP$ correspond to maximal flats in $\Es$.
This is still true in the non-Archimedean setting, 
i.e. replacing $\RR$ by an ultrametric valued field $\KK$,
in which case $\Es$ is a Euclidean building of type $A_{\N-1}$.

Configurations in projective spaces $\PP(\RR^\N)$ 
have been widely studied and used.
In particular, triples of flags in $\PP(\RR^\N)$ and 
their classical invariants 
(the triple ratio for $\N=3$), 
are the basic building blocks used by Fock and Goncharov 
to define generalized shearing coordinates
for higher Teichm{\"u}ller space, 
parametrizing positive representations 
of punctured surface groups in $\G=\SL_\N(\RR)$, 
see \cite{FoGoIHES}. 
But the geometric properties 
in the symmetric space or Euclidean building $\Es$ of these
configurations remain mysterious. 

In this article, we investigate the geometry 
of some simple ideal configurations
in a (not necessarily discrete) Euclidean building $\Es$ of type $A_2$, 
mainly the  generic triples of ideal  chambers,
and the relationship  with their projective geometry 
in the projective plane $\PP$.
Our first motivation is to use it to study actions 
of surface groups  on Euclidean buildings of type $A_2$,
and degenerations of Hitchin representations in $\SL_3(\RR)$ (see \cite{ParDegFG}).

The main result is a classification 
of ideal triples of chambers 
by the geometry of the five naturally associated flats in $\Es$,
in relation with  their triple ratio as triples of flags in $\PP$.
In the case where $\Es$ is a real tree 
(e.g. a Euclidean building of type $A_1$), 
any generic ideal triple
bounds a {\em tripod} in $\Es$, 
that is a convex subset consisting 
of union of three rays 
from a point $\xE\in\Es$ (the {\em center} of the tripod).
This is no longer the case in general
in higher rank buildings like $A_2$-buildings, 
and many types of configurations are possible.  
A special case was studied by A. Balser, 
who established a characterization 
of triples of points in $\bordinf\Es$ 
bounding a tripod in $\Es$ \cite{Balser08}, 
and used it to study convex rank $1$ subsets in $A_2$-buildings. 
We give here a complete and precise description.

We now get into more details.
Let $\Es$ be 
a real Euclidean building of (vectorial) type $A_2$, 
i.e. with model flat the Euclidean plane
\[\Aa = 
\left\{ \vA = (\vA_1,\vA_2,\vA_3) \in \RR^3 
  /\ \sum_\indN\vA_\indN=0\right\}\]
endowed with the finite reflection group $\W=\Sym_3$
acting by permutation of the coordinates.
Note that $\Es$ is not necessarily discrete (simplicial) 
nor locally compact, and possibly exotic.

The boundary at infinity of $\Es$ may be identified with 
the incidence graph of an associated projective plane 
$\PP=\PP_\infty(\Es)$, 
equipped with  
an $\RR$-valued additive cross ratio $\geombir$
(called a projective valuation in \cite{Tits86}) 
defined on quadruples of pairwise distinct collinear points in $\PP$
\cite{Tits86}.
In the algebraic case, 
i.e. when 
$\Es$ is the Bruhat-Tits building $\Es(\KK^3)$ 
associated with the group $\PGL(\KK^3)$ 
for some  ultrametric field $\KK$, 
the projective plane $\PP$ is $\PP(\KK^3)$ 
and $\geombir$ is the logarithm
\[\geombir=\log\abs{\Bir}\]
of the absolute value 
of the usual $\KK$-valued cross ratio $\Bir$ on $\PP(\KK^3)$, 
where conventions on cross ratios are taken such that
\[\Bir(\infty,-1,0, \Z)=\Z\]
in $\PP^1\KK=\KK\cup\{\infty\}$
(following \cite{FoGoIHES}).
We will then call $\geombir$ the {\em geometric} cross ratio 
and $\Bir$ the {\em algebraic} cross ratio to distinguish between them.

We now turn to ideal triples of chambers.
Let  $\Tau=(\F_1,\F_2,\F_3)$ 
be a triple of chambers at infinity of $\Es$. 
We denote by $\F_\ii=(\p_\ii,\D_\ii)$ the corresponding flag of $\PP$,
with $\p_\ii$ the point and $\D_\ii$ the line.
The set $\{1,2,3\}$ of indices 
will be canonically identified  with $\ZZ/3\ZZ$.
A triple   $\Tau=(\F_1,\F_2,\F_3)$ will be called {\em \genericTF} if
the flags $(\F_\ii)_\ii$ are pairwise opposite,
the points $(\p_\ii)_\ii$ are not collinear 
and the lines  $(\D_\ii)_\ii$ are not concurrent.

In the algebraic case $\PP=\PP(\KK^3)$ 
\genericTF{} triples of flags $\Tau=(\F_1,\F_2,\F_3)$
are classified 
by one $\KK$-valued invariant,  
the ({\em algebraic}) triple ratio 
(see for example \cite[\S 9.4]{FoGoIHES}),
that may be defined by:
\begin{equation}
  \label{eq:TriAsBir}
  \Tri(\F_1,\F_2,\F_3) =
  \Bir(\D_1, \p_1\p_2,\p_1\p_{23},\p_1\p_3)
\end{equation}
where $\p_{\ii\jj}=\D_\ii\cap \D_\jj$. 
We recall that it is invariant under cyclic permutations of 
$\Tau$,
and that reversing the order inverses the algebraic triple ratio: 
$\Tri(\ov{\Tau})=\Tri(\Tau)^{-1}$,
where $\ov{\Tau}=(\F_3,\F_2,\F_1)$.

In the general case, we introduce 
an invariant for \genericTF{} triples of flags in $\PP$,
analoguous to the algebraic triple ratio:
the {\em geometric triple ratio}, 
which still make sense when the building $\Es$ is exotic
(non algebraic), 
whereas the algebraic triple ratio is not defined anymore.
We define it as the  triple
\[\geomtri(\Tau)=(\geomtri_\indTriRap(\Tau))_{\indTriRap=1,2,3}\]
of the following cross ratios in $\PP$,
which are the cross ratios
obtained from the four lines $\D_1, \p_1\p_2, \p_1\p_{23}, \p_1\p_3$
by cyclic permutation of the three last one:
\[
\begin{array}{rl}
  \geomtri_1(\F_1,\F_2,\F_3) 
&=   \geombir(\D_1, \p_1\p_2, \p_1\p_{23}, \p_1\p_3)\\
  \geomtri_2(\F_1,\F_2,\F_3)
&=  \geombir(\D_1, \p_1\p_3, \p_1\p_2, \p_1\p_{23})\\
  \geomtri_3(\F_1,\F_2,\F_3)
&=   \geombir(\D_1, \p_1\p_{23}, \p_1\p_3, \p_1\p_2)
\end{array}
\;.\]

To simplify notations, we  denote from now on 
\[\gZ_\indTriRap=\geomtri_\indTriRap(\Tau)
  \text{ and } \gZ=(\gZ_1,\gZ_2,\gZ_3)=\geomtri(\Tau)\]
In the algebraic case,  we have $\PP=\PP(\KK^3)$ and 
the geometric triple ratio is obtained from the algebraic cross
ratio $\Z=\Tri(\Tau)\in\KK$ by:
\[  \begin{array}{rcl}
    \gZ_1&=&\log\abs{\Z}\\
    \gZ_2 &=&\log\abs{\frac{1}{1+\Z}}=-\log\abs{1+\Z}\\
    \gZ_3 &=& \log\abs{1+\Z^{-1}} \;.
    \end{array}
  \]
  The geometric triple ratio $\gZ$
enjoys the following properties. 
It is invariant  by cyclic permutations of the flags,
and
 changed to $(-\gZ_1,-\gZ_3,-\gZ_2)$ 
 by permutations reversing the cyclic order. 
We also have $\gZ_1+\gZ_2+\gZ_3=0$, and 
the stronger following property:
for all $\indTriRap\in\ZZ/3\ZZ$, 
if $\gZ_\indTriRap>0$ then $\gZ_{\indTriRap-1}=0$ and 
$\gZ_{\indTriRap+1}=-\gZ_\indTriRap<0$.
Note that the three natural cases: 
$\gZ\in\RR_+(0,1,-1)$, 
$\gZ\in\RR_+(-1,0,1)$, and
$\gZ\in\RR_+(1,-1,0)$
subdivide in two types, as the case $\gZ_1=0$ is invariant under
reversing the order of $\Tau$, 
whereas the two other  cases are exchanged.

We now turn to  the geometry inside the Euclidean building $\Es$. 
A \genericTF{} triple $\Tau=(\F_1,\F_2,\F_3)$ of ideal chambers  
defines five natural flats in $\Es$:
the three flats $\Aij=\App(\F_\ii,\F_\jj)$
containing the opposite chambers $\F_\ii$ and $\F_\jj$ in their boundaries,
the flat $\Ap=\App(\p_1,\p_2,\p_3)$ containing the triple
of ideal singular points $(\p_1,\p_2,\p_3)$ in its boundary,
and the similarly defined flat $\AD=\App(\D_1,\D_2,\D_3)$.
We will show that there are also six particular points in $\Es$ naturally
associated with the triple $\Tau$, that may be defined as the
orthogonal projections $\yTauE_\ii$ and $\yTauE_\ii^*$ (which happen
to be unique) of 
$\p_\ii$ and $\D_\ii$ on the flat $\Ajk$ 
where $\jj=\ii+1$ and $\kk=\ii+2$.
We say that $(\F_1,\F_2,\F_3)$ is of type {\em \typetripod{}}
if there exists a tripod in $\Es$  
joining the three (middle points of the) ideal 
chambers $(\F_1,\F_2,\F_3)$.  The set of centers of such
tripods is the intersection $I$ of the three flats $\Aij$.

We show that either 
the three flats $\Aij$ have nonempty intersection,
i.e. $(\F_1,\F_2,\F_3)$ is of type \typetripod,
or 
the two flats $\Ap$ and $\AD$ have non empty intersection $\Delta$,
which is then a {\em flat singular triangle} 
(that is, a triangle in $\Aa$ with singular sides)
(we then say that $(\F_1,\F_2,\F_3)$ is of type {\em  \typeflat{}}).
The two following results 
describe more precisely the two possible types,
and relate them with the points $\yTauE_\ii, \yTauE_\ii^*$ 
and the geometric triple ratio $\gZ$.
We denote by 
$\Cc=\{\vA\in\Aa/\ \vA_1 > \vA_2 >\vA_3\}$
the model Weyl chamber of $\Aa$ 
and we use the corresponding {\em simple roots coordinates} on $\Aa$,
that is $\vA=(\vA_1-\vA_2,\ \vA_2-\vA_3)$.

\begin{theorem}[Type \typetripod{}]
\label{theo- tripod}
The intersection $I=\App_{12} \cap \App_{23} \cap \App_{31}$ is
  nonempty if and only if 
$\gZ_1=   0$. 
Then $\gZ_2 \geq 0$ and 
there exist a unique pair $(\vstE,\vstE^*)$ in $\Es$ such that
\begin{myenum}
\item $\yTauE_1=\yTauE_2=\yTauE_3=\vstE$ 
and  $\yTauE_1^*=\yTauE_2^*=\yTauE_3^*=\vstE^*$ ;

\item $I$ is the segment $[\vstE,\vstE^*]$ ; 
\item $[\vstE,\vstE^*]$ is the unique shortest segment 
joining $\Ap$ to $\AD$.
\item 
\label{i- coords xx^*}
Identifying $\App_{\ii\jj}$ with $\Aa$ 
by a marked flat $\f:\Aa\mapsto \App_{\ii\jj}$ sending $\Cc$ to $\F_\jj$, 
in simple roots coordinates, we have
$\overrightarrow{\vstE\vstE^*}=(-\gZ_2,\gZ_2)$.
In particular $\vstE^*$ is on the ray $[\vstE,\p_{\ii\jj})$ 
from $\vstE$ to $\p_{\ii\jj}$. 
\end{myenum}
\end{theorem}

\begin{figure}[h]
\centering
\includegraphics[scale=0.8]{figures/Es_triple_type_tripod.fig}

  \begin{minipage}[b]{0.3\linewidth}
\centering
\includegraphics[scale=0.8]{figures/Es_Aij_tripod.fig}

\ In the flat $\Aij$.
  \end{minipage}
\hfill
  \begin{minipage}[b]{0.3\linewidth}
\centering
\includegraphics[scale=0.8]{figures/Es_Ap_tripod.fig}
    
\ In the flat $\Ap$.
  \end{minipage}
\hfill
  \begin{minipage}[b]{0.3\linewidth}
\centering
\includegraphics[scale=0.8]{figures/Es_AD_tripod.fig}

\ In the flat $\AD$.
  \end{minipage}

\caption{Type \typetripod}
\label{fig-  type tripod}
\end{figure}

\begin{theorem}[Type \typeflat{}]
\label{theo- triangle}
The intersection $\Ap \cap \AD$ is nonempty if and only if
($\gZ_2=0$ or $\gZ_3=0$), 
or, equivalently, if and only if  $\gZ_2\leq 0$. 
Then there exists a unique  flat singular triangle $\Delta\subset\Es$ with
  vertices $\vstE_1,\vstE_2,\vstE_3$
such that

\begin{myenum}
\item $\Ap \cap \AD = \Delta$.

\item
\label{i- Aij cap Aik chamber tip x_\ii}
    $\Aij\cap\Aik$ is the Weyl chamber from $\vstE_\ii$ to $\F_\ii$ ;

\item
\label{i- xi xj}
Let $\ii\in\{1,2,3\}$ and $\jj=\ii+1$.
In a marked flat $\f:\Aa\mapsto \Aij$ sending $\Cc$ to $\F_\jj$, 
in simple roots coordinates, we have
$\overrightarrow{{\vstE_\ii}{\vstE_\jj}}
=(\pospart{\gZ_1},\negpart{\gZ_1})$
where $\pospart{\gZ_1}=\max(\gZ_1,0)$ and $\negpart{\gZ_1}=\max(-\gZ_1,0)$.
In particular $\vstE_\jj$ is on the ray from
$\vstE_\ii$ to $\p_\jj$ (if $\gZ_1\geq 0$) 
or $\D_\jj$ (if $\gZ_1\geq 0$).

\item 
\label{i- Delta opp Fi en xi}
 The germs of Weyl chambers at $\vstE_\ii$ 
respectively defined  by  $\Delta$ and $\F_\ii$ 
are opposite (in the spherical building of directions at $\vstE_\ii$).
In particular there exists a flat containing
  $\Delta$, and containing $\F_\ii$ in its boundary.
\end{myenum}
Furthermore 
if $\gZ_1 \geq 0$ we have
 $\vstE_\ii=\yTauE_{\ii-1}=\yTauE^*_{\ii+1}$ for all $\ii$,
and if $\gZ_1 \leq 0$ we have $\vstE_\ii=\yTauE_{\ii+1}=\yTauE^*_{\ii-1}$ for
all $\ii$.
\end{theorem}

The intersections of each flat with the four other flats form a
partition (i.e. a covering with disjoint interiors),
which is  described in Figure
\ref{fig- type tripod} for the type \typetripod,  
and in Figure \ref{fig-  type triangle in positive case}
for the type  \typeflat{}
(see Proposition \ref{prop- in Aij
  - intersections}, Corollary \ref{coro- partition Aij} and
Proposition \ref{prop-
  in Ap - intersections}).

\begin{figure}[h]
\centering
\includegraphics[scale=0.8]{figures/Es_triple_type_triangle.fig}

  \begin{minipage}[b]{0.3\linewidth}
\centering
\includegraphics[scale=0.8]{figures/Es_Aij_triangle_pos.fig}

{In $\Aij$, with $\jj=\ii+1$.}
  \end{minipage}
\hfill
  \begin{minipage}[b]{0.3\linewidth}
\centering
\includegraphics[scale=0.8]{figures/Es_Ap_triangle_pos.fig}
    
In $\Ap$.
  \end{minipage}
\hfill
\begin{minipage}[b]{0.3\linewidth}
\centering
\includegraphics[scale=0.8]{figures/Es_AD_triangle_pos.fig}
    
\vskip -2.5ex
{\hskip 5em In $\AD$.}
  \end{minipage}

\caption{Type \typeflat, in the case where $\gZ_1  \geq 0$
(the case $\gZ_1\leq 0$ is obtained from the case $\gZ_1\geq 0$
by reversing the order of the flags $\F_\ii$, 
i.e. by exchanging $1$ and $3$ and $\ii$ and $\jj$ in the above pictures).
}
\label{fig-  type  triangle in positive case}
\end{figure}

\medskip

The special case where the hypotheses of both 
Theorems \ref{theo- tripod} and \ref{theo- triangle}
are satisfied corresponds to the case where  $\gZ_1=\gZ_2=\gZ_3=0$.
Then the five flats intersect in a unique point $\vstE$,
and, 
in the spherical building of directions at $\xE$,    
the triple of chambers induced by $\Tau=(\F_1,\F_2,\F_3)$ is \genericTF{}.
In particular we recover the characterization of \cite{Balser08} for
triples of points in $\bordinf\Es$ bounding a tripod in $\Es$. 
Note that M. Talbi  established some analogous  geometric classification
for interior triangles  in  discrete Euclidean buildings of type
$A_2$, see \cite{Tal06}.

Theorem \ref{theo- triangle} will be used in \cite{ParDegFG} to study
actions of punctured surface groups on Euclidean buildings of type
$A_2$. It allows us to give a metric interpretation, in the building,
of Fock-Goncharov parameters associated with  ideal triangulations.
We are then able to construct in $\Es$
an invariant weakly convex cocompact $2$-complex 
for large families of  actions. 
Theorem \ref{theo- triangle} enables 
us to associate to each triangle of the triangulation a flat singular
triangle in $\Es$, the complex is then obtained by 
connecting them gluing flat strips.
This allows to describe length spectra for large families of 
degenerations of convex projective structures
on surfaces.

We also show that generic quadruples of points
in $\PP$ (which will be called {\em projective frames})
define a nice center in $\Es$, with various characterizations, 
see Proposition \ref{prop- center of a projective frame}
(this result generalizes to higer rank $\RR$-buildings of type
$A_{\N-1}$).

\paragraph{{\bfseries Aknowledgments}}
I would like to thank Fr{\'e}d{\'e}ric Paulin 
for useful discussions and comments. 
I also want to thank the members 
of the  Institut Fourier for their support.

\section{Preliminaries}
\label{s- prelim}

\subsection{The model flat 
\texorpdfstring{$(\Aa,\W)$}{(A,W)}
of type 
\texorpdfstring{$A_{\N-1}$}{A\_N} 
}
Let $\N\geq 2$ be an integer.  
The {\em model flat} of type $A_{\N-1}$ is  the vector space 
$\Aa=\RR^\N/\RR(1,\ldots,1)$, endowed with the action of
the {\em Weyl group}
 $\W=\Sym_\N$ acting on $\Aa$ by permutation of
coordinates (finite reflection group). 
We denote by $\classA{\vRN}$
the projection in $\Aa$ of a vector $\vRN$ in $\RR^\N$.
The vector space $\Aa$ may be identified with the hyperplane 
$\{\vA=(\vA_1,\ldots,\vA_\N)\in\RR^\N/\ \sum_\indN
\vA_\indN=0\}$ of $\RR^\N$.
Recall that a vector in $\Aa$ is called {\em singular} if it belongs
to one the hyperplanes $\vA_\indN=\vA_\indNp$, and {\em regular} otherwise. 
A {\em (open) (vectorial) Weyl chamber} of $\Aa$ is a connected component of
regular vectors.
We will call a {\em sector} a more general convex cone in $\Aa$, in
particular the closed convex cone formed by the union of the closed
Weyl chambers containing a given singular ray.
The {\em model Weyl chamber} 
is the simplicial cone  
\[\Cc=\{\vA\in\Aa/\ \vA_1 > \cdots >\vA_\N\}\;.\]
Its closure $\Cb$ is a strict fundamental domain for the action of
$\W$ on $\Aa$.
Recall that two nonzero vectors $\vA$ and $\vA'$ of $\Aa$ 
are called {\em opposite} if $\vA'=-\vA$.
Similarly, two Weyl chambers $\Ch$ and $\Ch'$ of $\Aa$ 
are {\em opposite} if $\Ch'=-\Ch$.
The {\em type}  of a vector $\vA\in\Aa$ is its 
projection  (modulo
$\W$) in $\Cb$.

We denote by $\bord \Aa$ the sphere of unitary vectors in $\Aa$,
identified with the set $\PP^+(\Aa)=(\Aa-\{0\})/\RR_{>0}$ 
of rays issued from $0$,
and by $\bord: \Aa-\{0\} \to \bord \Aa$ the corresponding projection.
The {\em type (of direction)}  of 
a nonzero vector $\vA \in \Aa$ 
is its canonical projection 
in $\bord \Cb$.

We denote by
$(\eRRNo_1,\ldots,\eRRNo_\N)$ the canonical basis of $\RR^\N$.
\label{s- def type of vertices in bord Aa}
For $d=1,\ldots,\N-1$, we will say that  a nonzero vector in $\Aa$ (or a
point in the sphere $\bord\Aa$) is {\em singular of type $d$} if
its canonical projection in $\bord\Cb$ is
$\classA{\eRRNo_1+\cdots+\eRRNo_d}$.

The  {\em simple roots} (associated with $\Cc$) are the following linear forms on $\Aa$
\[\rac{\indN} : \vA\mapsto \vA_\indN-\vA_{\indN+1}\]
for $\indN=1,\ldots,\N-1$.
The set of simple roots is denoted by $\RacSimples$.
We will also use the root 
$\rac\N : \vA \mapsto\vA_\N-\vA_1 $
satisfying
\[\rac1+\cdots+\rac\N =0\;.\]

\begin{figure}[h]
\centering
  \includegraphics[scale=0.75]{figures/Aa.fig} 
  \caption{The model flat $\Aa$ of type
    $A_2$ (for $\N=3$), and simple roots coordinates.
    The arrows denote the singular directions of type $1$.}
  \label{fig- Aa}
\end{figure}

The vector space $\Aa$ is endowed with 
the unique $\W$-invariant Euclidean scalar product, 
which is  well defined up to homothety (induced by the
standard Euclidean scalar product of $\RR^\N$).
We will normalize it 
by requiring that the simple roots have unit norm, 
i.e. the distance between the two hyperplanes 
with equation $\rac\indN=0$ and $\rac\indN=1$
is $1$ for one (all) $\indN$. 
When $\dim\Aa=1$, 
we will  identify $\Aa$ 
with $\RR$ by the basis $\{\classA{\eRRNo_1}\}$, i.e. 
by the map from  $\sR \mapsto  \classAcoordsll{\sR}{0}$
from $\RR$ to $\Aa$, which is an
isometry in the above normalization.

\subsection{Projective spaces}
\label{ss- projective spaces} 
We here collect the notations and vocabulary for projective spaces,
which will be used throughout this article.
We refer to \cite[\S 6.2]{Tits74}.
Let $\PP$ be a projective space  of dimension $\N-1$, with $\N\geq 2$.
We denote by $\AllFlags(\PP)$ the set of flags of $\PP$,
that is  increasing sequences $(\V_1,\ldots,\V_{\M})$ 
of proper linear subspaces of $\PP$.
We denote by $\PP^*$ the dual projective space, whose set of points is the set  of hyperplanes of $\PP$.

Two maximal flags $(\V_1,\ldots,\V_{\N-1})$ and $(\V'_1,\ldots,\V'_{\N-1})$
are {\em opposite} if they are in generic position, that is if
$\V_i\oplus \V'_{n-i}=\PP$ for all $i$.
A finite subset $\p_1, \ldots, \p_\M$ in $\PP$,  with $2\leq \M \leq \N$, 
is called {\em \independent{}} if it is not contained in any linear subspace of
dimension $\M-2$ of $\PP$. 
Then it is contained in a unique $(\M-1)$-dimensional linear subspace of
$\PP$, which will be denoted by $\p_1\oplus \cdots \oplus \p_{\M}$.
When $\M=2$, we will also denote the line $\p\oplus\q$ by $\p\q$.

A {\em frame}
of $\PP$ 
is a \independent{} $\N$-tuple. 
A {\em projective frame} in $\PP$
is a $(\N+1)$-tuple $(\p_0,\p_1, \ldots,\p_\N)$  
of points in $\PP$
in generic position, 
i.e. such that  
the induced $\N$-tuple
$(\p_0,\ldots, \widehat{\p_\indN},\ldots, \p_\N)$ is a frame in $\PP$ for all $\indN$.

If $\p$ is a point in $\PP$, 
we denote by $\PP/\p$
the  set  of lines through $\p$,
which is a projective space of dimension $\N-2$ whose linear subspaces
are the linear subspaces of $\PP$ containing $\p$.
The {\em projection at $\p$} is the corresponding
projection $\projP_\p:\q \mapsto \p\q$
from $\PP-\{\p\}$ to $\PP/\p$.
If $\p$ is a point of $\PP$ and $H\subset \PP$ an hyperplane with $\p\notin H$,
 then the projection $\projP_\p$ induces a canonical 
isomorphism  $\projP_{H\p}:H\isomto \PP/\p$ (called {\em perspectivity}).

Note that if $\Frame=(\p_1,\ldots, \p_\M)$ is \independent{} in $\PP$, 
then its projection  
$\projP_{\p_1}(\Frame)=(\p_1\p_2,\ldots, \p_1\p_\M)$  
at $\p_1$
is \independent{} in  $\PP/\p_1$.
In particular the projection of a (projective) frame at
one of its points is a (projective) frame.

\subsection{Spherical buildings of type 
\texorpdfstring{$A_{\N-1}$}{A\_{N-1}}
 and associated projective spaces}
\label{ss- proj space ass to a spherical An building}
See \cite[\S 6]{Tits74}.
A spherical building $\B$ of type $A_{\N-1}$ 
is the building of flags 
of an associated projective space $\PP=\PP(\B)$ of dimension $\N-1$. 
For $d=0,1,\ldots,\N-1$, 
the set of linear subspaces of dimension $d$ of $\PP$ 
identifies with the subset of vertices of type $d+1$
 of $\B$.
In particular, the projective space  $\PP$ itself
is identified with the set of vertices of type $1$ of $\B$,
and the dual projective space $\PP^*$ is identified with
the set of vertices of type $\N-1$.

In the algebraic case, 
that is when $\B$ is the spherical building of flags 
of
some vector space $\V$ of dimension $\N$ over a field $\KK$, 
then $\PP=\PV$.

A basic fact is that 
frames in $\PP$ correspond to apartments of $\B$.

Recall that, in (the geometric realization
modeled on $(\bord\Aa,\W)$ of) 
a spherical building, any two points (resp. chambers)
are contained in a common apartment,
and that they are  {\em opposite} if they
are opposite in that apartment, 
that is, 
for two points $\xi$ and $\xi'$, 
if and only if $\sphericalangle(\xi,\xi')=\pi$ 
for the canonical metric $\sphericalangle$ on $\B$.
Note that $\p\in\PP$ and $H\in\PP^*$ 
are opposite if and only if 
$\sphericalangle(\p,H)=\pi$, 
if and only if $\p\notin H$.
Two chambers are opposite if and only if they are opposite
as maximal flags in $\PP$.
In particular, in the type $A_2$ case, two chambers 
$\F_1=(\p_1,\D_1)$, $\F_2=(\p_2,\D_2)$ are opposite 
if and only if $\p_1\notin \D_2$ and $\p_2 \notin \D_1$.

For any simplex $\sigma$ of $\B$ the {\em residue}  $\Star(\sigma)$ 
of $\sigma$ is  the
spherical building formed by 
the simplices of $\B$ containing $\sigma$. 
If $H$ is a hyperplane of $\PP$, 
the residue $\Star(H)$ of $H$ in $\B$ 
is the subset of  flags of $\PP$ containing $H$.
It canonically identifies 
with the spherical building $\AllFlags(H)$ of flags of $H$ 
by the map $(\V_1,\ldots,\V_{\M},H)\mapsto (\V_1,\ldots,\V_{\M})$.
The residue $\Star(\p)$ of a point $\p$ in $\PP$
identifies canonically
with the flag building $\AllFlags(\PP/\p)$ of $\PP/\p$ 
by the map 
$(\V_1=\p,\ldots,\V_{\M})\mapsto (\V_2/\p,\ldots,\V_{\M}/\p)$.
If $\p\notin H$ then the projection $\projP_\p$ induces a canonical 
isomorphism  $\projP_{H\p}:\Star(H)\isomto
\Star(\p)$ of spherical buildings (perspectivity).

\subsection{Euclidean buildings}
Euclidean buildings considered in this article are
 (not necessarily discrete) Euclidean buildings of type
$A_{\N-1}$.
We refer for example to \cite{ParImm} for the definition and
properties of Euclidean buildings we use below 
(see also \cite{Tits86}, \cite{KlLe97}, \cite{Rousseau09}). 
Recall that a {\em Euclidean building  of type
$A_{\N-1}$} is a $\CAT(0)$ metric space 
$\Es$ endowed with a (maximal)
collection $\mathcal{A}$ of
isometric embeddings $\f:\Aa \to \Es$ called  {\em marked apartments},
or {\em marked flats} by analogy with  Riemannian symmetric spaces,
satisfying the following properties
\begin{mydescription}
  \item[{\bfseries (A1)}] $\mathcal{A}$ is invariant by precomposition by $\Waff$;
  \item[{\bfseries (A2)}] If $\f$ and $\f'$ are two marked flats, then the
    transition map $\f^{-1}\circ \f'$ is  the restriction of an element of
 $\Waff$;
  \item[{\bfseries (A3')}] Any two rays of $\Es$ are initially 
 contained in a common marked flat;
\end{mydescription}
where $\Waff$ denotes the subgroup of all affine isomorphisms of $\Aa$
with linear part in $\W$.  The {\em flats} 
 and the {\em Weyl chambers}) 
of $\Es$ are the images by the marked flats of $\Aa$ and  $\Cc$, respectively.

\subsubsection*{Algebraic case}
\label{s- Es(V)}
Let $\KK$ be an ultrametric field, i.e. a field endowed with 
an ultrametric absolute value $\abs{\cdot}$ (not necessarily
discrete). 
When $\V$ is a finite $\N$-dimensional vector space over $\KK$, 
we denote by  $\Es=\Es(\V)$ the Euclidean building associated with
$\G=\PGL(\V)$.
 We refer for example to \cite{ParImm} 
for the model of norms for $\Es$ (see \cite{GoIw63}, \cite{BrTi84}).
To each  basis $\eVb$ of $\V$ is then
associated a marked flat $\f_\eVb:\Aa\to \App_\eVb\subset\Es$,
such that, 
if $\aG$ is an element of $\G$ with diagonal
matrix $\diag(\aK_1,\ldots,\aK_\N)$ in the basis $\eVb$, 
then $\aG$ translates the flat $\App_\eVb$  by the  vector
\[\nu(\aG)=\classA{(\log\abs{\aK_\indN})_\indN}\]
in $\Aa$ (identifying the flat $\App_\eVb$ with the model flat
$\Aa$ through the marking $\f_\eVb$).

\medskip

{\em From now to Section \ref{s- A-Busemann cocycle},
$\Es$ will denote a Euclidean building of type $A_{\N-1}$.}

\subsection{Spherical building and projective space at infinity}
The CAT(0) boundary $\bordinf\Es$ of $\Es$ is the geometric
realization modeled on $(\bord\Aa,\W)$ of a spherical building of type
$A_{\N-1}$ whose chambers are the boundaries of the Weyl
chambers of $\Es$, and whose apartments are the boundaries of the flats of $\Es$.
It will be identified with the building of flags on the associated
projective space $\PP=\PP_\infty(\Es)$,
whose points are the vertices of type $1$ of $\bordinf \Es$. 
If $\chp$ and $\chm$ are opposite ideal chambers, then
we denote by 
$\App(\chm,\chp)$ the unique flat {\em joining} $\chm$ to $\chp$ in $\Es$,
that is, containing $\chm$ and $\chp$ in its boundary.
If $\Frame$ is a frame of $\PP$ or $\PP^*$, then 
there is a unique flat  $\App(\Frame)$ of $\Es$ containing $\Frame$ in its boundary.

\subsection{Local spherical building and projective space at a point}
Recall that, in Euclidean buildings, 
two (unit speed) geodesic segments issued from
a common point $\xE$ have zero angle 
if and only if 
they have same germ at $\xE$ (i.e. coincide in a neighborhood of $\xE$).
A {\em direction at $\xE\in\Es$} is 
a germ of nontrivial geodesic segment from $\xE$.
A direction, geodesic segment, ray or line 
has a well-defined {\em type (of direction)} in $\bord\Cb$, 
which is its canonical
projection (through a marked flat) in $\bord\Cb$.
It is called {\em singular} or {\em regular} accordingly.

The {\em space of directions} 
at  $\xE$ of $\Es$ 
is the quotient space of non trivial geodesic segments from $\xE$ for
this relation, with the induced angular metric,
and is denoted by  $\TangS_\xE \Es$. 
We denote by   $\TangS_\xE:\Es-\{\xE\}\to \TangS_\xE \Es$,
$\yE \to \TangS_\xE\yE$, 
the associated projection.
Its extension to the boundary at infinity 
will also be denoted by $\TangS_\xE: \bordinf\Es \to \TangS_\xE \Es$,
$\xi \to \TangS_\xE\xi$
and called the {\em canonical projection}.

The space of directions $\TangS_\xE \Es$ inherits the structure 
of a spherical $A_{\N-1}$-building,
whose apartments are the germs $\TangS_\xE \App$ at $\xE$
of the flats $\App$ of $\Es$ passing through $\xE$,
and whose chambers are the germs $\TangS_\xE \Ch$ at $\xE$ 
of the Weyl chambers $\Ch$ of $\Es$ with vertex $\xE$ 
(see for example \cite{ParImm}).
The canonical projection $\TangS_\xE: \bordinf\Es \to \TangS_\xE \Es$
sends chambers to chambers 
(and, more generally, simplices to simplices)
and  
preserves the type of points. 

The {\em local projective space} $\PP_\xE = \PP_\xE(\Es)$ {\em at} $\xE$ 
is the projective space of dimension $\N-1$ associated with 
 the spherical building $\TangS_\xE \Es$ of type $A_{\N-1}$
(see \S \ref{ss- proj space ass to a spherical An building}).
Its underlying set is the set  of vertices of type $1$ of $\TangS_\xE\Es$.

The canonical projection $\TangS_\xE: \bordinf\Es \to \TangS_\xE \Es$
induces  (by restriction to vertices) 
a surjective morphism (of projective spaces)  
$\TangS_\xE:\PP \to \PP_\xE$
from the projective space at infinity $\PP$ 
to the local projective space $\PP_\xE$ at $\xE$.
Note that, in particular, if $\Frame$ is a frame of $\PP$, then 
\label{s- x is in App(Frame) iff TangS_xFrame is a frame}
$\xE$ belongs to the associated flat $\App(\Frame)$ 
if and only if 
$\TangS_\xE(\Frame)$ is a frame of $\PP_\xE$.

\subsection{Transverse spaces at infinity}
\label{s- transverse spaces at infinity}
See for example \cite[\S 8]{Tits86}, 
\cite[1.2.3]{Leeb00}, \cite[\S 4]{MSVM14-I}.
Let $\xi$ be a vertex of $\bordinf\Es$ of type $1$ or $\N-1$,
i.e. either a point $\p$ in the projective plane at infinity $\PP$ or
a hyperplane $H$ of $\PP$. 

The {\em transverse space} $\Es_\xi$ at $\xi$
may be defined, from the metric viewpoint 
(as in \cite[1.2.3]{Leeb00}),
as 
the quotient space of the set of all rays to $\xi$  by the
pseudodistance $\dxi$ given by
\[\dxi(\ray_1,\ray_2)=\inf_{\tR_1,\tR_2}\deucl( \ray_1(\tR_1),
\ray_2(\tR_2) ) \;.\]%
We denote by  $\projE_\xi:\Es\to\Es_\xi$ the canonical projection 
(which maps $\xE$ to the class of the unique ray from $\xE$ to $\xi$).
The space $\Es_\xi$  is a Euclidean building 
of type $A_{\N-2}$,
whose flats are the projections to $\Es_\xi$ of the flats
of $\Es$ containing a ray to $\xi$.
In particular, when $\Es$ is of type $A_2$,
the transverse space $\Es_\xi$  is an $\RR$-tree, and we will call it the
{\em transverse tree} at $\xi$.

In the algebraic case, i.e. when $\Es=\Es(\V)$, 
the transverse space $\Es_H$ canonically
identifies with the building $\Es(H)$ of $H$, 
where $H$ is seen as  an hyperplane of $\V$,
and $\Es_\p$ identifies with $\Es(\V/\p)$, 
where $\p$ is seen as a $1$-dimensional subspace of $\V$.

The spherical building $\bordinf \Es_\xi$ at infinity  of $\Es_\xi$ 
identifies canonically with the residue $\Star(\xi)$ of $\xi$.
In particular,  if  $\p$ is a point in $\PP$, 
the projective space at infinity of $\Es_\p$ 
identifies with  $\PP/\p$, 
and if $H$ is an hyperplane of $\PP$, 
the projective space at infinity of $\Es_H$ identifies with $H$.

\label{projection of a flat and frames}
If $\Frame=(\p_1, \ldots, \p_\N)$  is a frame in $\PP\subset\bordinfEs$, 
then the projection on $\Es_{\p_1}$ of the flat
$\App(\p_1, \ldots, \p_\N)$ 
is the flat defined by the projection
$\projP_{\p_1}(\Frame)=(\p_1\p_2,\ldots, \p_1\p_\N)$ of the frame
$\Frame$, i.e. $\projE_{\p_1}(\App(\Frame))=\App(\projP_{\p_1}(\Frame))$.

We now describe the canonical isomorphism 
$\projE_\ximp:\Es_\xim \isomto \Es_\xip$ 
for  opposite points $\xim$ and $\xip$ of $\bordinfEs$.
The union  $\Esximp$ of all geodesics joining $\xim$ to $\xim$ 
is a convex closed subspace and a subbuilding,
whose flats are the flats
of $\Es$ containing a  geodesic joining $\xim$ to $\xim$
(see \cite[prop. 4.8.1]{KlLe97} and \cite[2.2.1]{ParComp}).
We denote by $\Esximp=\Es^\ximp\times \RR$
the canonical decomposition (see \cite[1.2.10]{ParRepcr}).
The restriction of the projection $\projE_\xip$ to $\Esximp$ is surjective
and factorizes through the projection on the first factor, 
inducing a canonical isomorphism of Euclidean buildings
$\Es^\ximp \isomto \Es_\xip$.
We similarly have a isomorphism
$\Es^\ximp \isomto \Es_\xim$,
so it induces a canonical isomorphism 
$\projE_\ximp:\Es_\xim \isomto \Es_\xip$.
\label{s-prel- persectivities comes from an 
isomorphism of the transverse trees}%
It is easy to see that the map $\projE_\ximp$ 
extends to the boundaries at infinity of 
$\Es_\xim$ and $\Es_\xip$
by the canonical isomorphism of spherical buildings 
$\projP_\ximp:\Star(\xim) \isomto \Star(\xip)$ (perspectivity).

\subsection{The 
\texorpdfstring{$\Aa$}{A}-valued 
Busemann cocycle}
\label{s- A-Busemann cocycle}
Let $\ch$ be a chamber at infinity of $\Es$. 
We now  define the  $\Aa$-valued {\em Busemann cocycle}
\[\Buso_\ch:\Es\times\Es\to\Aa\]
associated to $\ch$.
It can be simply defined from canonical
retractions as 
\[\Bus{\ch}{\xE}{\yE}:=\ret(\yE)-\ret(\xE)\]
where $\ret : \Es\to \Aa$  is any canonical retraction centered at
$\ch$, sending $\ch$ to $\bord\Cc$ (see \cite[Prop. 1.19]{ParImm})).
More precisely,
the Buseman cocycle at $\ch$ is characterized by the property: 
\[\Bus{\ch}{\f(\vA)}{\f'(\vA')}=\vA'-\vA\]
for any two marked flats $\f,\f':\Aa\to \Es$
sending $\bord\Cc$ to $\ch$ and such that $\f=\f'$ on some subchamber
of $\Cc$.

We clearly
have
\[\Bus{\ch}{\xE}{\zE}=\Bus{\ch}{\xE}{\yE}+\Bus{\ch}{\yE}{\zE}\;.\]
When $\dim\Aa=1$, it coincides with the usual Busemann cocycle, 
\label{ss- usual Bus cocycle}%
which is defined for $\xi \in \bordinf\Es$ by
\[\Bus{\xi}{\xE}{\yE}=\lim_{\zE\ra \xi}\deucl(\xE,\zE)-\deucl(\yE,\zE)\;.\]
In the type $A_2$ case, the simple root coordinates of
$\Aa$-valued Busemann cocycles
 may be determined by projecting in transverse trees at infinity, 
using the following relations (using the normalization of the metric).
\begin{equation}
  \label{eq- projections and Busemann cocycle}
  \begin{array}{rl}
    \rac1(\Bus{(\p,\D)}{\xE}{\yE})
    & = \Bus{\p} {\projE_\D(\xE)} {\projE_\D(\yE)}\\
    \rac2(\Bus{(\p,\D)}{\xE}{\yE})
    & = \Bus{\D} {\projE_\p(\xE)} {\projE_\p(\yE)} \;.
  \end{array}
\end{equation}

\bigskip

We now turn to cross ratios.

\subsection{Cross ratio on the boundary of a tree}
\label{ss- geometric cross ratio on the boundary of a tree}
See \cite[\S 7]{Tits86}, 
and in a more general setting \cite{Otal92}, \cite{Bourdon96}.
In this section, we suppose that $\Es$ is 
an $\RR$-tree.
Given three distinct ideal points $\xi_1,\xi_2,\xi_3$ in $\bordinfEs$,
we denote by $\centre(\xi_1,\xi_2,\xi_3)$ 
the {\em center} of the ideal triple $(\xi_1,\xi_2,\xi_3)$, 
that is the unique common intersection point 
of the three geodesic lines joining two of the three points.
Note that $\centre(\xi_1,\xi_2,\xi_3)$ is the (orthogonal)
projection of $\xi_3$ on the geodesic joining $\xi_1$ to $\xi_2$.
We denote by $\Bus{\xi}{\xE}{\yE}$ the Busemann cocycle 
(see \S \ref{ss- usual Bus cocycle}).
\noindent {\advance\linewidth by-5cm
\begin{minipage}{\linewidth}
Define the {\em cross ratio} of four pairwise distinct 
points $\xi_1$, $\xi_2$, $\xi_3$, $\xi_4$ in  $\bordinfEs$ 
by
\[ \geombir(\xi_1, \xi_2, \xi_3,\xi_4)
=\frac{1}{2}(\ell_{12}-\ell_{23}+\ell_{34}-\ell_{41})\]
where $\ell_{\ii\jj}$ is the length of the geodesic in $\Es$ 
from $\xi_\ii$ to $\xi_\jj$ after removing disjoint
fixed horoballs centered at each $\xi_\kk$. It does not
dependend on the choice of the horoballs since the horoballs
centered at a given point are  equidistant along the rays to that point.
\end{minipage}
\hfill
 \begin{minipage}{4cm}
\centering\includegraphics[scale=0.75]{figures/Tree_geom_crossratio.fig}
\end{minipage}
}

The cross ratio naturally extends to 
 {\em \nondegenerateQ{}} quadruples,
that are 
quadruples $(\xi_1,\xi_2,\xi_3,\xi_4)$ 
{\em without triple point}  
(i.e. any three of the points are not equal),
which is equivalent to the following condition:
\begin{equation}
\label{eq- nondegenerateQ}
  (\xi_1 \neq \xi_4 \mbox{ and } \xi_2 \neq \xi_3)
  \mbox{ or }
  (\xi_1 \neq \xi_2 \mbox{ and }  \xi_3 \neq \xi_4)
\;.
\end{equation}
  
We then set
\[
\geombir(\xi_1,\xi_2,\xi_3,\xi_4) 
= \left\{
\begin{array}{cl}
0       & \mbox{ when } \xi_1 = \xi_3 \mbox{ or } \xi_2 = \xi_4\\
\minfty & \mbox{ when } \xi_1 = \xi_2 \mbox{ or } \xi_3 = \xi_4\\
\pinfty & \mbox{ when } \xi_1 = \xi_4 \mbox{ or } \xi_2 = \xi_3
\end{array}
\right.%
\;.
\]
We now recall some basic properties that we will use.

\noindent {\advance\linewidth by-5cm
\begin{minipage}{5cm}
\centering\includegraphics{figures/Tree_geom_crossratio_centers_tripods.fig}
\end{minipage}
\hfill
\begin{minipage}{\linewidth}
The cross ratio may be read inside the tree on the oriented geodesic
from $\xi_3$ to $\xi_1$, as the oriented distance 
$\ovra{\xT\yT}$ 
from the center $\xT$  
of the ideal triple
$(\xi_3,\xi_1,\xi_2)$ to the center $\yT$
of the ideal triple $(\xi_3,\xi_1,\xi_4)$:
\begin{equation}
\label{eq- geombir and centers of tripods}
  \geombir(\xi_1, \xi_2, \xi_3,\xi_4)
  =\ovra{\xT\yT}=\Bus{\xi_1}{\xT}{\yT}
\;.
\end{equation}
\end{minipage}
}

\medskip

The cocycle identity 
is
\[\geombir(\xi_1, \xi_2, \xi_3,\xi_4)+\geombir(\xi_1, \xi_4,\xi_3,\xi_5)=
\geombir(\xi_1, \xi_2, \xi_3,\xi_5)\;.\]

The cross ratio $\geombir$ is left unchanged by the double
transpositions
\label{ss- geombir (12)(34)}  
and changed to $-\geombir$ by the transpositions
$(13)$ and $(24)$.
We now consider the behaviour under cyclic permutations of the three last terms. 
We have
\begin{equation}
  \label{eq- sum of 3-cyclic perm of geom bir is 0} 
    \geombir(\xi_1, \xi_2, \xi_3, \xi_4)
  + \geombir(\xi_1, \xi_4, \xi_2, \xi_3)
  + \geombir(\xi_1, \xi_3, \xi_4, \xi_2) = 0
  \;.
\end{equation}
Moreover, the following     
{\em ultrametricity} property (specific to the case of trees) 
is easy to prove using \eqref{eq- geombir and centers of
  tripods} (see \cite[\S 7, prop. 3]{Tits86}):
\begin{equation}
\label{eq- geombir is ultram}
\begin{array}{ll}
\mbox{If } \geombir(\xi_1, \xi_2, \xi_3,\xi_4)>0, 
&\mbox{then } \geombir(\xi_1, \xi_3, \xi_4,\xi_2)=0\\
&\mbox{and }
\geombir(\xi_1, \xi_4, \xi_2,\xi_3)=-\geombir(\xi_1,
\xi_2,\xi_3,\xi_4)
\end{array}
\;.
\end{equation}
Note that \eqref{eq- geombir is ultram}
is equivalent 
(under \eqref{eq- sum of 3-cyclic perm of geom bir is  0}) 
to
\begin{equation}
\label{eq- geombir is ultram via max}
\geombir(\xi_1, \xi_2, \xi_3,\xi_4) 
\leq \max(0,-\geombir(\xi_1, \xi_4,  \xi_2,\xi_3))%
\;.\end{equation}
which in the algebraic case 
follows from the symmetry properties of the cross ratio 
under $3$-cyclic permutations \eqref{eq- sym bir 1342}.

\subsection{Algebraic case: link with usual cross ratio}
Suppose that $\Es$ is the tree 
$\Es(\V)$ associated with  a 2-dimensional
vector space $\V$ over an ultrametric field $\KK$ 
(see Section \ref{s- Es(V)}).
Then $\bordinfEs$ identifies with the projective line $\PV$.

The usual cross ratio $\Bir$ on $\PV$ of a 
\nondegenerateQ{} quadruple of points 
(see \eqref{eq- nondegenerateQ}) 
is defined by (following the convention of \cite{FoGoSPC},
and taking values in  $\KK\cup\{\infty\}$)
\begin{equation}
\label{eq- def cross ratio}  
\Bir(\aK_1,\aK_2,\aK_3,\aK_4)
=\frac{(\aK_1 -\aK_2 )(\aK_3 -\aK_4 )}
      {(\aK_1 -\aK_4 )(\aK_2 -\aK_3 )}
\end{equation}
in any affine chart $\PV\isomto \KK\cup\{\infty\}$,  
so that $\Bir(\infty,-1,0,\aK)=\aK$.

The  cross ratio $\geombir$ 
defined in Section \ref{ss- geometric cross ratio on the boundary of a tree}
will then be called the {\em geometric} cross ratio, to distinguish it
from $\Bir$, which will be called the  {\em  algebraic} cross ratio. 
They are  then related  as follows:
\begin{equation}
  \label{eq- link geom and alg cross ratio} 
 \geombir(\xi_1, \xi_2, \xi_3,\xi_4)
=\log\abs{\Bir(\xi_1, \xi_2,  \xi_3,\xi_4)}
\;.\end{equation}
\begin{proof}
Let $\xT_4=\centre(\xi_3,\xi_1,\xi_2)$ and
$\xT_2=\centre(\xi_3,\xi_1,\xi_4)$.
In a suitable basis $\eVb=(\eV_1,\eV_2)$ of $\V$, 
we have in homogeneous coordinates
  $\xi_1=[1:0]$, $\xi_3=[0:1]$, $\xi_2=[-1:1]$
 and $\xi_4=[\bK:1]$, 
where $\bK=\Bir(\xi_1, \xi_2,\xi_3,\xi_4)$. 
Then $\g=
  \begin{pmatrix}
-\bK&0\\
0&1
  \end{pmatrix}$
fixes $\xi_1$ and $\xi_3$ and sends $\xi_2$ to $\xi_4$.
Hence $\g(\xT_4)=\xT_2$. In the flat $\App(\xi_3,\xi_1)$ identified with
$\Aa=\RR^2/\RR(1,1)$ by the marked flat $\f_\eVb$,
 we have
$\ovra{\xT_4\xT_2}=\nu(\g)
=\classA{(\log\abs{\bK},0)}$, 
hence $\ovra{\xT_4\xT_2}=\log\abs{\bK}$ as needed.
\end{proof}

We recall that the algebraic cross ratio $\Bir$
satisfies the following symmetry properties:
It is left unchanged by the double transpositions
and changed to $\Bir^{-1}$ by  the transpositions $(13)$ and $(24)$. 
Furthermore we have 
an additional symmetry under $3$-cycles 
not satisfied by the geometric cross ratio:
\begin{equation}
\label{eq- sym bir 1342}
    \begin{array}{rl}
      \Bir(\aK_1,\aK_3,\aK_4,\aK_2)&=-1-\Bir(\aK_1,\aK_2,\aK_3,\aK_4)^{-1}\\
      \Bir(\aK_1,\aK_4,\aK_2,\aK_3)&=-(1+\Bir(\aK_1,\aK_2,\aK_3,\aK_4))^{-1}  
    \end{array}
\;.\end{equation}

\subsection{Cross ratio on the boundary 
  of an \texorpdfstring{$A_2$}{A\_2}-Euclidean building}
\label{ss- geometric cross ratio on a projective plane}
See \cite{Tits86}. 
Let $\Es$ be a Euclidean building of type $A_2$, and $\PP$ the
associated projective plane at infinity. 

Let $(\p_1,\p_2,\p_3,\p_4)$ be a \nondegenerateQ{} quadruple
of points of $\PP$ on a common line $\D$. 
Then their {\em cross ratio} $\geombir(\p_1,\p_2,\p_3,\p_4)$
(i.e. {\em projective valuation} in \cite{Tits86})
is by definition 
their cross ratio as ideal points of the transverse tree $\Es_\D$.
The cross ratio of a \nondegenerateQ{} quadruple of lines in $\PP$ 
passing through a common point $\p$ 
is similarly defined 
as their cross ratio 
as ideal points of the transverse tree $\Es_\p$.

The main additional property is that perspectivities preserve cross
ratio, which follows from the fact that perspectivities extend
isometries between the transverse trees (see Section \ref{s-prel-
  persectivities comes from an isomorphism of the transverse trees}):

\begin{proposition}
Let $\p$ be a point of $\PP$ 
and $\D$  a line of $\PP$ 
with $\p \notin \D$.
The canonical isomorphisms (perspectivities)
$\projP_{\p\D}: \Star(\D)\isomto \Star(\p)$, $\q \mapsto \p\q$ 
and 
$\projP_{\D\p}: \Star(\p)\isomto \Star(\D)$, $L \mapsto \D\cap L$,
preserve the cross ratio $\geombir$, i.e.
\begin{equation*}
\geombir(\p_1, \p_2,\p_3, \p_4)=\geombir(\p\p_1, \p\p_2,\p\p_3, \p\p_4)  
\end{equation*}
\begin{equation*}
\geombir(\D_1, \D_2,\D_3, \D_4)
=\geombir(\D\cap\D_1,\D\cap\D_2,\D\cap\D_3, \D\cap\D_4)  
\end{equation*}
\vskip -5ex
\qed\end{proposition}

\section{Some basic ideal configurations}

\subsection{
  Extension of orthogonal projection 
  to the boundary 
  in CAT(0) spaces
}
In this section $\Es$ is a general  CAT(0) metric space, and we prove the following basic property:
the usual orthogonal projection onto a proper convex subset $\Y\subset \Es$ 
 extends to the boundary outside 
the closed $\frac{\pi}{2}$-neighborhood of $\bordinf\Y$ for the Tits metric
(note that the projection is no longer unique).
This property is quite elementary 
but we did not see it 
in the classical litterature, 
so we include the proof.
We refer to the  book \cite{BrHa99} for CAT(0) spaces.

We denote by $\bordinf\Es$ the CAT(0) boundary of $\Es$, and by
$\sphericalangle_{Tits}(\xi,\eta)$ the Tits angle between two ideal points $\xi,\eta\in\bordinf\Es$.
For a subset $A$ of $\bordinf\Es$,
we define $\sphericalangle_{Tits}(\xi,A)=\inf_{\eta\in A}\sphericalangle_{Tits}(\xi,\eta)$.
\begin{definition}
  Let $\Y$ be a subspace of $\Es$ and $\xi\in\bordinf\Es$ an ideal point. 
We say that a point  $\xE\in \Y$ is an {\em orthogonal projection of $\xi$ on $\Y$} if
  $\sphericalangle_\xE (\xi, \yE) \geq \frac{\pi}{2}$ for all $\yE\in \Y-\{\xE\}$.
\end{definition}

\begin{proposition}
  \label{prop- CAT0 - projection of ideal points}
  Let $\Y$ be a convex subspace of a CAT(0) space
  $\Es$ which is proper for the induced metric, 
  and $\xi$ in $\bordinf\Es$. 
  Suppose that $\sphericalangle_{Tits}(\xi,\bordinf\Y)> \frac{\pi}{2}$. 
  Then there exists an orthogonal projection $\xE$ of $\xi$ on $\Y$.
\end{proposition}

\begin{proof}
  Consider a sequence $(\xE_\n)$ converging to $\xi$ in $\Es$, and let $\yE_\n$ be the
  orthogonal projection of $\xE_\n$ on $\Y$.
  If $(\yE_\n)_{\n\in\NN}$ is not bounded then, up to passing to a
  subsequence, 
  $\yE_\n$ converges to $\eta$ in $\bordinf\Y$. 
  Then for any fixed $\yE$ in $\Y$ we have  
  $\sphericalangle_\yE(\xi,\yE_\n) \leq \frac{\pi}{2}$ for all $\n$,
  hence $\sphericalangle_\yE(\xi,\eta) \leq \frac{\pi}{2}$.
  Therefore  $\sphericalangle_{Tits}(\xi,\eta) \leq \frac{\pi}{2}$.
  Thus $(\yE_\n)_{\n\in\NN}$ is bounded, hence, since   $\Y$ is proper, 
  it has a converging subsequence, and the limit
  point $\xE$ is then an orthogonal projection of $\xi$ on $\Y$.
\end{proof}

\subsection{Centers of generic 
\texorpdfstring{$(\N+1)$}{(\N+1)}-tuples.}
In this section, we show that the notion of center of ideal triples
in trees extends to Euclidean buildings of type $A_{\N-1}$, for
generic $(\N+1)$-tuples of points (or hyperplanes) in the associated
projective space at infinity (Proposition \ref{prop- center of a projective frame}).

Let $\Es$ be a Euclidean building of type $A_{\N-1}$, 
and  $\PP$ be its projective space at infinity  
(i.e., the set of singular points of type $\stypep$ in $\bordinf\Es$,
see Section \ref{s- prelim}).  
Recall from Section \ref{ss- projective spaces} that 
a {\em projective frame} 
in a projective space of dimension $\N-1$ 
is  a generic $(\N+1)$-tuple of points.

We first observe that 
the orthogonal projection of a point of $\PP$  on a flat of $\Es$ 
exists under a simple  necessary and sufficient condition.

\begin{proposition}
\label{prop- proj p sur A}
Let $\App$ be a flat of $\Es$ 
and $\p\in\PP$. 
Let $(\p_1,\ldots,\p_\N) = (\bordinf\App) \cap \PP$ 
be the points of type $\stypep$ in $\bordinf\App$.
Then $p$ admits an orthogonal projection on $\App$ 
if and only if
$(\p,\p_1,\ldots,\p_\N)$ is a projective frame.
\end{proposition}

The analoguous property is also valid for points $H\in\PP^*$.
Note that these properties also hold in
symmetric spaces of type $A_{\N-1}$.

\begin{proof}
  If $\p\in H$
  for some hyperplane $H$ in $\PP^*\cap\bordinf\App$, then
  $\p$ and $H$ are in a common chamber of the spherical building
  $\bordinf\Es$, and,
  as the diameter $d$ of the model spherical Weyl chamber $\bord \Cb$ 
  is stricly less that $\pi/2$ (for the angle metric), 
  we have
  $\sphericalangle_{Tits}(\p,H)<{\pi/2}$, hence the orthogonal projection
  does not exist.
  Else, for every hyperplane $H$ in $\PP^*\cap\bordinf\App$, we have
  $\p\notin H$, hence $\sphericalangle_{Tits}(\p,H)=\pi$,
  which implies that   since $\sphericalangle_{Tits}(\p, \eta) \geq  \pi-d > \pi/2$ 
  for all $\eta\in\bordinf\App$,
  and the orthogonal projection exists by
  Proposition \ref{prop- CAT0 - projection of ideal  points}.
\end{proof}

We now turn to the main result of this section.

\begin{proposition}
  \label{prop- center of a projective frame}  
  Let $\Frame=(\p_0, \p_1,\ldots, \p_\N)$  be a projective frame 
  in $\PP\subset\bordinfEs$.
  For each $\indN\in\{0,\ldots, \N\}$ let $\App_\indN$ be the unique 
  flat of $\Es$ through
  $(\p_0,\ldots, \widehat{\p_\indN},\ldots, \p_\N)$.
  There exists a unique point $\xE\in\Es$
  satisfying the following equivalent conditions.

  \begin{myenum}

  \item 
    \label{i- x in intersection of all flats}
    $\xE \in \cap_\indN \App_\indN$ ;

  \item 
    \label{i-  angles at x are pi}
    For all $\indN$ and for all  $H$ in $\bordinf\App_\indN\cap\PP^*$  
    the angle $\sphericalangle_\xE (\p_\indN,H)$ is $\pi$ ;
  \item 
    \label{i- directions at x}
    The $(\N+1)$-tuple   $\TangS_\xE\Frame=(\TangS_\xE
    \p_\indN)_{\indN=0,\ldots,\N}$ 
    of directions at $\xE$
    form a projective frame in $\PP_\xE$ ;

  \item 
    \label{i- for all ii, x is the proj on App_ii}
    For all $\indN$, the point  $\xE$ is an orthogonal projection 
    of $\p_\indN$ on the flat  $\App_\indN$ ;

  \item 
    \label{i- for one ii, x is the proj on App_ii}
    There exists $\indN$ such that $\xE$ is an orthogonal  
    projection of $\p_\indN$ on  $\App_\indN$.
  \end{myenum}

  We will call $\xE$  
  the {\em center}  of the projective frame $\Frame=(\p_0, \p_1,\ldots, \p_\N) $ 
  and denote it by $\centre(\p_0, \p_1,\ldots, \p_\N)$ or $\centre(\Frame)$.
\end{proposition}

\begin{figure}[h]
\centering
  \includegraphics[scale=1]{figures/Es_ideal_simplex.fig}
\caption{
The center $\xE\in\Es$
of a projective frame  $(\p_1,\p_2,\p_3$, $\p_4)$ (for $\N=3$).
}
\end{figure}

\begin{proof}
  The existence of $\xE$, as an orthogonal  projection of $\p_0$ on  $\App_0$,
  is ensured by  Proposition \ref{prop- proj p sur A}.
  
  For $\indN\neq\indNp$, denote by $H_{\indN\indNp}$ the hyperplane
  $\oplus_{\indNpp\neq\indN,\indNp}\; \p_\indNpp$ in the projective space $\PP$.
  Let $\xE\in\Es$. 
  Conditions \ref{i- directions at x} 
  and \ref{i- x in intersection of all  flats} 
  are equivalent 
  (see Section \ref{s- x is in App(Frame) iff TangS_xFrame is a frame}).

  We first show 
  \ref{i- x in intersection of all flats}
  $\Rightarrow$ 
  \ref{i-  angles  at x are pi}:
  Fix $\indN$ and  $H\in\PP^*$  in $\bordinf\App_\indN$.
  The opposite of $H$ in $\bordinf\App_\indN$ is some $\p_\indNp$. 
  Then $H=H_{\indN\indNp}$, 
  so $H$ is also the opposite of $\p_\indN$
  in the apartment $\bordinf\App_\indNp$. As $\xE\in\App_\indNp$, we then have
  $\sphericalangle_\xE (\p_\indN,H)=\pi$. 
  We now prove 
  \ref{i-  angles  at x are pi} 
  $\Rightarrow$
  \ref{i- directions at x}: 
  First recall that 
  for $\p\in\PP$ and $H\in\PP^*$, we have 
  $\sphericalangle_\xE (\p_\indN,H)=\pi$ if and only if
  $\TangS_\xE\p\notin\TangS_\xE H$ in the projective space $\PP_\xE$.
  So \ref{i-  angles  at x are pi} means that 
  $\TangS_\xE\p_\indN \notin \TangS_\xE H_{\indN\indNp}$ for all $\indN\neq\indNp$.
  Let $\U_\indN$ be the minimal linear subspace 
  of the projective space $\PP_\xE$ 
  containing $\TangS_\xE\p_0,\ldots,\TangS_\xE\p_\indN$. 
  Then, for $\indN\leq \N-1$, 
  we have  that $\TangS_\xE\p_\indN$ is not in $\U_{\indN-1}$, 
  else $\TangS_\xE\p_\indN$ would belong to $\TangS_\xE H_{\indN,\indN+1}$. 
  Hence
  $(\TangS_\xE\p_0,\ldots,\TangS_\xE\p_\indN)$ is \independent{} in $\PP_\xE$ 
  by induction on $\indN$.
  Therefore $(\TangS_\xE\p_0,\ldots,\TangS_\xE\p_{\N-1})$ is a frame, 
  and \ref{i- directions at x} follows by permuting the $\p_\indN$.

  We now prove 
  \ref{i-  angles at x are pi} 
  $\Rightarrow$
  \ref{i- for all ii, x is the proj on  App_ii}.
  Let $\indN\in\{0,\ldots,\N\}$.
  Let $\vE\in \TangS_\xE\App_\indN$.
  Let $\Ch\subset\App_\indN$ be a closed Weyl chamber with vertex $\xE$ 
  containing $\vE$. 
  Let $H\in\PP^*$ be the singular point of type $\stypeH$ in $\bordinf\Ch$. 
  Then $\sphericalangle_\xE (\p_\indN,H)=\pi$, 
  hence $\sphericalangle_\xE (\p_\indN, \vE)\geq \pi- d >\frac{\pi}{2}$, 
  as the diameter $d$ of  $\bord \Cb$ 
  is stricly less that $\pi/2$.

  \ref{i- for all ii, x is the proj on  App_ii}
  $\Rightarrow$ 
  \ref{i- for one ii, x is the proj on  App_ii}
  is clear.
  Assume now that \ref{i- for one ii, x is the proj on App_ii} holds.
  For $\indNp\neq \indN$ in $\{0,\ldots,\N\}$,
  as $\sphericalangle_\xE(\p_\indN,H_{\indN\indNp})\geq \frac{\pi}{2}$,
  the direction $\TangS_\xE\p_\indN$ is not in a
  closed chamber of $\TangS_\xE\Es$ containing $\TangS_\xE H_{\indN\indNp}$.
  Hence by type considerations we must have $\sphericalangle_\xE(\p_\indN,H_{\indN\indNp})=\pi$.
  So \ref{i-  angles at x are pi} holds.

  So the equivalence of all assertions is proven.
  We now prove the uniqueness of $\xE$. 
  Suppose that $\xEp$ is another point of $\Es$ with the same
  properties, and $\xEp \neq \xE$. We proved above that we have then
  $\sphericalangle_\xE (\p_\indN, \xEp)> \frac{\pi}{2}$
  and $\sphericalangle_{\xEp} (\p_\indN, \xE)> \frac{\pi}{2}$, which is impossible.
\end{proof}

We now state some properties of centers of projective frames.
Consider   
 a projective frame $\Frame=(\p_0, \p_1,\ldots,\p_\N)$   in $\PP$, 
and let $\xE\in\Es$ be its center.
Let $\App_\indN=\App(\p_0,\ldots, \widehat{\p_\indN},\ldots, \p_\N)$
be the $\N+1$ associated flats in $\Es$.
We first describe the intersections of the flats $\App_\indN$ with $\App_0$.

\begin{proposition}
  \label{prop- Ai cap Aj and partition Ai}
For $\indN=1\ldots\N$, 
let $\SectE_\indN$ be
the sector with base-point $\xE$ 
on 
$\{\p_1,\ldots, \widehat{\p_\indN},\ldots, \p_\N\}$,
i.e.
  the convex hull of the rays from $\xE$ to these points.
And let  
$H_\indN=\p_1 \oplus \cdots \oplus \widehat{\p_\indN} \oplus\cdots \oplus \p_\N$
denote the point in $\bordinf\App_0$ opposite to $\p_\indN$.
For $\indN \in \{1,\ldots,\N\}$, we have:

\begin{myenum}

\item 
\label{i- dir p_0}
Let $\yE$ be an interior point of $\SectE_\indN$.
Then  $\TangS_\yE \p_0 = \TangS_\yE \p_\indN$.

\item 
\label{i- p_0 opp H_i}
For  $\yE\in\App_0$, we have
$\yE\in\App_0\cap\App_\indN$ 
if and only if    
$\TangS_\yE\p_0$ is opposite to $\TangS_\yE H_\indN$;

\item 
$\App_0\cap\App_\indN= \SectE_\indN$ ;
\end{myenum}
In particular, the intersections $\App_0\cap\App_\indN$, $\indN=1\ldots\N$,   
form a partition (i.e. a covering with disjoint interiors) of
$\App_\indN$.
\end{proposition}

Note that the sector $\SectE_\indN$ is
the union of the Weyl chambers of the flat $\App_0$ based at $\xE$
and containing the singular  ray to $H_\indN$.

\begin{proof}
The  inclusion $\SectE_\indN\subset \App_0\cap \App_\indN$ is clear 
since $\xE\in  \App_0\cap\App_\indN$ 
and $\p_\indNp$ is in 
$\bordinf\App_0  \cap \bordinf \App_\indN$ 
for $\indNp\neq \indN$ in $\{1,\ldots,\N\}$.

If  $\yE$ is an interior point of $\SectE_\indN$, then
in the local spherical building $\TangS_\yE\Es$ at $\yE$,
we have that 
$\TangS_\yE\p_0 \in \TangS_\yE\App_0$.
Moreover, $\yE\in \App_\indN$ 
as previously observed, 
so $\TangS_\yE\p_0$ is opposite to  $\TangS_\yE H_\indN$  
(in $\TangS_\yE\App_\indN$).
Hence $\TangS_\yE\p_0$ is equal to 
the opposite of $\TangS_\yE H_\indN$ in $\TangS_\yE\App_0$, which is 
$\TangS_\yE\p_\indN$, proving \ref{i- dir p_0}.

We now prove \ref{i- p_0 opp H_i}:
In
$\PP_\yE$, 
the points 
$(\TangS_\yE\p_1,\ldots, \TangS_\yE\p_\N)$ 
form a  frame
(since $\yE\in\App_0$).
Hence the $\N-1$ points
$(\TangS_\yE\p_1,\ldots,\widehat{\TangS_\yE\p_\indN},\ldots, \TangS_\yE\p_\N)$ 
are \independent{}
Therefore  
$(\TangS_\yE\p_0,\ldots,\widehat{\TangS_\yE\p_\indN},\ldots, \TangS_\yE\p_\N)$ 
is a frame in $\PP_\yE$
(i.e. $\yE\in\App_\indN$)
if and only if
$\TangS_\yE\p_0 \notin \TangS_\yE H_\indN$.

We finish by proving the remaining inclusion 
$\App_0\cap \App_\indN\subset \SectE_\indN$:
The $\SectE_\indN$ clearly form a partition of $\App_0$.
So it is enough to
prove that
 that  $\App_0\cap \App_\indN$ does not meet
the interior of $\SectE_\indNp$ for $\indNp\neq \indN$.
Else, at such a point $\yE$, by \ref{i- dir p_0}, 
we would have $\TangS_\yE \p_0=\TangS_\yE \p_\indNp$,
which  is not opposite to $\TangS_\yE H_\indN$, 
providing a contradiction.
\end{proof}

The following proposition shows that 
the notion of center of projective frames
behaves well with respect to projections to transverse spaces at
infinity.

\begin{proposition}
  \label{prop- proj of centers to transverse spaces} 
  For each $\indN$, the projection of $x$ in the transverse building at
  infinity $\Es_{\p_\indN}$ is 
  the center of the projective frame of
  $\bordinf\Es_{\p_\indN}$ formed by the projections
  $\projP_{\p_\indN}(\p_\indNp)=\p_\indN\p_\indNp$ of the $\p_\indNp$, $\indNp\neq \indN$,
  that is:
  \[ \projE_{\p_\indN}(\centre(\p_0, \p_1, \ldots, \p_\N))
  =\centre
  (\p_\ii\p_0, \p_\ii\p_1,\ldots, \widehat{\p_\ii\p_\ii},\ldots,\p_\ii\p_\N)
  \;.\]
\end{proposition}

\begin{proof}
  For all  $\indNp\neq \indN$,
  the ray from $\xE$ to $\p_\indN$ 
  is in the flat  $\App_\indNp$ hence
  its projection $\projE_{\p_\indN}(\xE)$ 
  in the transverse building  $\Es_{\p_\indN}$
  is in $\projE_{\p_\indN}(\App_\indNp)$, 
  which is the flat defined by the frame
  $\projP_{\p_\indN}(\p_\indNpp)=\p_\indN\p_\indNpp$, $\indNpp\neq\indN,\indNp$.
\end{proof}

In the algebraic case, 
i.e. when  $\Es$ is the Euclidean buiding $\Es(\V)$ 
associated with some vector space $\V$ of dimension $\N$
 over an ultrametric field $\KK$, 
we have the following characterization of the center 
as a norm on $\V$.

\begin{proposition}
  Let $\Frame=(\p_0, \p_1,\ldots, \p_\N)$  be a projective frame in $\PP=\PV$.
  The center  of $\Frame$
  is the norm $\normV$ on $\V$ canonically associated 
  to any basis $\eVb=(\eV_\indN)_{\indN=1,\ldots \N}$ of $\V$ such that
  $\p_\indN=\classPV{\eV_\indN}$ for $1\leq \indN \leq \N$ and
  $\p_{0}=\classPV{\eV_1+\cdots+\eV_\N}$ in $\PV$, 
  i.e. the norm defined by 
  \[\normV(\sum_{\indN=1}^\N\aK_\indN\eV_\indN)=\max_{1\leq \indN \leq \N}\abs{\aK_\indN}\;.\]
\end{proposition}

\begin{proof}
  Let $\eVb=(\eV_1,\ldots,\eV_\N)$ be a basis of $\V$ such that
  $\p_\indN=\classPV{\eV_\indN}$ and
  $\p_{0}=\classPV{\eV_1+\cdots+\eV_\N}$ in $\PV$.
  Let $\normV$ 
  be the associated canonical norm on $\V$.
  We clearly have $\normV \in \App_{0}$ by the definition of marked
  flats in the model of norms.
  Let $\g$ be the element of
  $\GL(\V)$ sending the basis $\eVb$ to the basis
  $(\eV_1,\ldots,\eV_{\N-1},\eV_1+\cdots+\eV_\N)$.
  Then $\g$ 
  preserves the norm $\normV$ 
  and sends $\App_0$ to $\App_\N$ and 
  hence $\normV$ is in the flat $\App_\N$.
  Permuting the basis $\eV$, we similarly get  that $\normV$ is in 
  the flat $\App_\indN$ for all $\indN\neq 0$.
\end{proof}

\begin{remark}
  By duality, the similar properties hold for generic $(\N+1)$-tuples
  (projective frames)
  in  $\PP^*\subset\bordinf\Es$.
\end{remark}

\subsection{Projecting  two ideal points onto a flat}
From now on we return to the case where $\N=3$ (type $A_2$).

\begin{proposition}
\label{prop- 1 plat, 2 points}
Let $(\p_1,\p_2,\p_3)$ be a \independent{} triple in $\PP$.
Let $\p,\q$ be two points in $\PP$, 
in generic position relatively to the $\p_\ii$ 
(i.e. not on any of the lines $\p_\ii\p_\jj$).
Denote by $\xE$ and $\yE$ the respective orthogonal projections of $\p$ and $\q$ on
the flat $\App=\App(\p_1,\p_2,\p_3)$.
Identify $\App$ with $\Aa$ by a marked flat
sending $\bord\Cc$ to $(\p_1,\p_1\p_2)$.
Then the roots coordinates of $\ovra{\xE \yE}$ are given by the three
natural cross ratios at the vertices of the triangle:
\[\rac1( \ovra{\xE \yE})
=  \geombir(\p_3\p_1, \p_3\p, \p_3\p_2 , \p_3\q),\]
\[\rac2( \ovra{\xE \yE})
= \geombir(\p_1\p_2, \p_1\p,  \p_1\p_3 , \p_1\q ),\]
\[\rac3( \ovra{\xE \yE})
= \geombir(\p_2\p_3, \p_2\p,  \p_2\p_1 , \p_2\q )\;.\]
\end{proposition}

The analogous dual result holds for projections of two lines of $\PP$ on a flat
(exchanging
the roles of points and lines in $\PP$).

\begin{proof}
Projecting on the transverse tree $\Es_{\p_1}$ in direction $\p_1$,
 we have
\[\rac2( \ovra{\xE \yE})
=\rac2( \Bus{(\p_1,\p_1\p_2)} {\xE} {\yE} )
=\Bus{\p_1\p_2} {\projE_{\p_1}(\xE)} {\projE_{\p_1}(\yE)}\]
by  
\eqref{eq- projections and Busemann cocycle}.
Since the projections of $\xE$ and $\yE$ on the tree $\Es_{\p_1}$ are the
respective centers of the ideal triples $(\p_1\p_2, \p_1\p_3, \p_1\p)$
and $(\p_1\p_2, \p_1\p_3, \p_1\q)$ 
(Proposition \ref{prop- proj of centers to transverse spaces}),
we have 
\[\Bus{\p_1\p_2} {\projE_{\p_1}(\xE)} {\projE_{\p_1}(\yE)}
=  \geombir(\p_1\p_2, \p_1\p,  \p_1\p_3 , \p_1\q )\]
by \eqref{eq- geombir and centers of tripods},
hence
$\rac2( \ovra{\xE \yE})
= \geombir(\p_1\p_2, \p_1\p,  \p_1\p_3 , \p_1\q )$.
The remaining assertions follow by applying cyclic permutation, since 
\[\rac1(\Bus{(\p_1,\p_1\p_2)}{\xE}{\yE})
=\rac2(\Bus{(\p_3,\p_3\p_1)}{\xE}{\yE})\]
\[\rac3(\Bus{(\p_1,\p_1\p_2)}{\xE}{\yE})
=\rac2(\Bus{(\p_2,\p_2\p_3)}{\xE}{\yE})\; .\]
\vskip -3ex
\end{proof}

For the projections of a point and a line, we have the following result.

\begin{proposition}
\label{prop- 1 flat p D}
  Let $\F_-=(\p_-,\D_-)$ and $\F_+=(\p_+,\D_+)$ be two oppposite flags
  in $\PP$
and $\App$ the flat in $\Es$ joining them, identified with $\Aa$
by a marked flat sending $\bord\Cc$ to $\F_+$.
Let $\p$ be a point and $\D$ a line 
in $\PP$ in generic position with respect to $\F_-$ and $\F_+$,
(i.e.
$\p$ does not belong to  any of the lines $\p_-\p_+$, $\D_-$,
$\D_+$, and
$\D$ does not contain any of the points $\D_-\cap\D_+$, $\p_-$,
$\p_+$).

Denote by $\xE$ and $\xE^*$ the respective orthogonal projections of $\p$ and $\D$ 
on $\App$. Then in simple roots coordinates we have 
\[\ovra{\xE \xE^*}=(\gZ_-,\gZ_+),\]
\begin{align*}
\mbox{with }
\gZ_-
&=  \geombir(\p_+, \D_+\cap (\p_-\p),  \D_+\cap\D_- , \D_+\cap \D)\\
&=  \geombir(\D_-, \p_- \oplus (\D_+\cap \D),  \p_-\p_+ , \p_-\p )\\
\mbox{and }
\gZ_+
&=  \geombir(\p_-, \D_-\cap \D,  \D_-\cap\D_+ , \D_-\cap (\p_+\p))\\
&=  \geombir(\D_+, \p_+\p,  \p_+\p_- , \p_+ \oplus (\D_-\cap \D) )
\;.\end{align*}
\end{proposition}

\begin{figure}[h]
  \centering
  \includegraphics[scale=0.75]{figures/Es_proj_p_D_sur_A.fig}

  \caption{Projecting a point and a line  on a flat.}
\label{fig-projPointLineOnFlat}  
\end{figure}

\begin{proof}
  See Figure \ref{fig-projPointLineOnFlat}.
  The projection of $\xE$ on the transverse tree $\Es_{\p_-}$ 
is  the center of the ideal triple 
$(\p_-\p_+, \p_-(\D_-\cap\D_+), \p_-\p)$, and
the projection of $\xE^*$ on the tree $\Es_{\D_+}$ 
is  the center of the ideal triple 
$(\p_+, \D_+\cap\D_-,\D_+\cap\D)$ 
(Proposition \ref{prop- proj of centers to transverse spaces}).
As $\xE$ lies on a geodesic from $\p_-$ to $\D_+$, 
we have
\begin{align*}
\projE_{\D_+}(\xE) 
&= \projE_{\D_+,\p_-} ( \projE_{\p_-}(\xE) )\\
&= \projE_{\D_+,\p_-} ( \centre(\p_-\p_+, \p_-(\D_-\cap\D_+), \p_-\p) )\\
&=\centre (\p_+, \D_-\cap\D_+, \D_+ \cap (\p_-\p) )\;.
\end{align*}
Then projecting on the transverse tree $\Es_{\D_+}$ we have
\[\rac1( \ovra{\xE \xE^*})
= \Bus{\p_+}{\projE_{\D_+}(\xE)} {\projE_{\D_+}(\xE^*)}
=  \geombir(\p_+, \D_+\cap (\p_-\p),  \D_+\cap\D_- , \D_+\cap \D)\]
as needed. The remaining assertions have identical proofs.
\end{proof}

\section{Triple ratio of a triple of ideal chambers%
}
\label{s- triple ratio}

In this section, we introduce  
the {\em (geometric) triple ratio} of a \nondegenerateTF{} 
triple of ideal chambers in a
a real Euclidean building $\Es$ of type $A_2$,
establish its basic properties, 
and the links with the
usual $\KK$-valued (algebraic) triple ratio of triples of flags
(see e.g. \cite{FoGoSPC})
in the algebraic case $\PP=\PP(\KK^3)$.

We first give a precise definition
of {\em \nondegenerateTF{}} and {\em \genericTF} triples of flags 
in an arbitrary projective plane  $\PP$.

\subsection{\NondegenerateTF{} and \genericTF{}
  triples of flags}
\label{s- nondegenerateTF}
Let $\PP$ be a projective plane
and $\Tau=(\F_1,\F_2,\F_3)$ be a triple of flags
 $\F_\ii=(\p_\ii,\D_\ii)$ in $\PP$.
We will  denote by $\p_{\ii\jj}$ the point $\D_\ii\cap \D_\jj$
 (resp. $\D_{\ii\jj}$ the line $\p_\ii\p_\jj$), when defined.

The natural nondegeneracy condition on the triple $(\F_1,\F_2,\F_3)$ 
for the triple ratios  to be well defined is the following:
\begin{center}
(\ND)
either $\forall \ii,\ \p_\ii \notin \D_{\ii+1}$ 
or $\forall \ii,\ \p_\ii \notin \D_{\ii-1}$. 
\end{center}
This condition is clearly equivalent to:
the points are pairwise distinct, the lines are pairwise distinct,
none of the points is on the three lines 
(i.e. $\D_\ii\cap \D_\jj \neq \p_\kk$ for all $\{\ii,\jj,\kk\}=\{1,2,3\}$)
and none of the lines contains the three points
(i.e. $\p_\ii \p_\jj \neq \D_\kk$ for all $\ii,\jj,\kk$). 
We will then say that the triple $(\F_1,\F_2,\F_3)$ is 
{\em \nondegenerateTF}.

It is easy to check that the triple $\Tau$ defines then 
a \nondegenerateQ{} quadruple $(\D_\ii, \p_\ii\p_\jj, \p_\ii\p_{\jj\kk},
\p_\ii\p_\kk)$ of lines through each point $\p_\ii$,
and a \nondegenerateQ{} quadruple
$(\p_\ii, \D_\ii\cap\D_\jj, \D_\ii\cap\D_{\jj\kk}, \D_\ii\cap\D_\kk)$ 
of points  on each line $\D_\ii$.
The  triple of flags  $\Tau=(\F_1,\F_2,\F_3)$ is {\em \genericTF} if
the flags $\F_\ii=(\p_\ii,\D_\ii)$ are pairwise opposite, 
the points $(\p_\ii)_\ii$ are not collinear 
and the lines  $(\D_\ii)_\ii$ are not concurrent.
In particular, $\Tau$ is then \nondegenerateTF{}, and the induced
quadruples of points on each line (resp. of lines through each point)
are generic (i.e. pairwise distinct).

\subsection{Algebraic triple ratio}

When $\PP=\PP(\KK^3)$ is the projective plane associated 
with an arbitrary field $\KK$,
 the algebraic triple ratio of a \nondegenerateTF{} 
triple of flags $\Tau=(\F_1,\F_2,\F_3)$ 
(see Section \ref{s- nondegenerateTF}), 
with values in $\KK\cup\{\infty\}$,  is defined by
(see  \cite[\S 9.4]{FoGoIHES})
\[
\Tri(\F_1,\F_2,\F_3)
=
\frac{\Dt_1(\pt_2)\Dt_2(\pt_3)\Dt_3(\pt_1)}{\Dt_1(\pt_3)\Dt_2(\pt_1)\Dt_3(\pt_2)}
\]
where $\pt_\ii$ is any vector in $\KK^3$ representing $\p_\ii$ and
$\Dt_\ii$ is any linear form in $(\KK^3)^*$ representing $\D_\ii$, and
$\F_\ii=(\p_\ii,\D_\ii)$.
It is invariant under cyclic permutation of the flags 
and inversed by reversing the order
\[\Tri( \F_3, \F_2, \F_1 )=\Tri(\F_1,\F_2,\F_3)^{-1} \;.\]
It may be expressed as the following cross ratio
\begin{align}
\label{eq- triple ratio as a cross ratio}
\Tri(\F_1,\F_2,\F_3)
&=\Bir(\D_1, \p_1 \p_2, \p_1 \p_{23}, \p_1 \p_3)\;.
\end{align}

\subsection{Geometric triple ratio}

We suppose now that the projective plane $\PP$ is the projective plane
at infinity of some a real Euclidean building $\Es$ of type
$A_2$, possibly exotic. Let $\geombir$ be the associated geometric cross ratio on $\PP$ 
(see Section \ref{ss- geometric cross ratio on a projective plane}).
Let   $\Tau=(\F_1,\F_2,\F_3)$ be a \nondegenerateTF{} 
triple of ideal chambers of $\Es$,
i.e. a \nondegenerateTF{}  triple of flags $\F_\ii=(\p_\ii,\D_\ii)$ in $\PP$. 

The idea is to define the geometric triple ratio 
of  $\Tau$ by analogy with the expression of the algebraic triple ratio as a cross ratio
\eqref{eq- triple ratio as a cross ratio}, replacing $\Bir$ by
$\geombir$,
in such a way that, in the algebraic case, the geometric
triple ratio of a triple $\Tau$ with algebraic triple ratio $\Z$ should
be $\log\abs{\Z}$.
But for the purpose of geometric classification, this geometric cross
ratio $\geombir(\D_1, \p_1\p_2, \p_1\p_{23}, \p_1\p_3)$ alone will not
retain enough information.
In particular, in contrast to the algebraic
cross ratio, it does not determine the geometric cross ratios obtained
from the original $4$-tuple  by cyclic permutations of the three last
arguments, which in the algebraic case are  
$\log\abs{1+\Z^{-1}}$ and $-\log\abs{1+\Z}$,
see \eqref{eq- sym bir 1342},
and have geometric significance.
For example, in the algebraic case, it will not distinguish between two
triples $\Tau$ and $\Tau'$ with respective algebraic triple ratios $\Z=-1$
and $Z'=-1+a$ with $\abs{a}<1$. 

In order to retain this information we
define the  {\em geometric triple ratio} of $\Tau$ 
as the triple
\[\geomtri(\Tau)=(\geomtri_\indTriRap(\Tau))_{\indTriRap=1,2,3}\]
where
\[
\begin{array}{rcl}
  \geomtri_1(\F_1,\F_2,\F_3) &=& 
  \geombir(\D_1, \p_1\p_2, \p_1\p_{23}, \p_1\p_3)\\
  \geomtri_2(\F_1,\F_2,\F_3)&=&
  \geombir(\D_1, \p_1\p_3, \p_1\p_2, \p_1\p_{23})\\
  \geomtri_3(\F_1,\F_2,\F_3)&=& 
  \geombir(\D_1, \p_1\p_{23}, \p_1\p_3, \p_1\p_2)
\end{array}
\]
are the geometric cross ratios obtained from
$(\D_1, \p_1\p_2,\p_1\p_{23}, \p_1\p_3)$ 
by cyclic permutation of the three last lines.
Note these cross ratios are well defined, since the four lines
$\D_1$, $\p_1\p_2$, $\p_1\p_{23}$, $\p_1\p_3$ 
are well defined and form a \nondegenerateQ{}
quadruple of lines through $\p_1$ 
(see Section \ref{s- nondegenerateTF} above).
The following proposition gathers the properties of the geometric
triple ratio, and show in particular that this invariant is in fact
$1$-dimensional, as it takes values in one of the three rays
$\RR_+(0,1,-1)$,  $\RR_+(-1,0,1)$, and $\RR_+(1,-1,0)$.

\begin{proposition}
  \label{prop- props of geomtri}
The following hold.
  \begin{myenum}
  \item 
    \label{i- geomtri: invariant by cyclic perm}
    The geometric triple ratio is invariant by cyclic permutations of
    the flags, i.e.
    for $\indTriRap=1,2,3$, 
    \[\geomtri_\indTriRap(\F_2,\F_3,\F_1)
    =\geomtri_\indTriRap(\F_1,\F_2,\F_3) \; ;\]

  \item Exchanging two flags, we have
    \label{i- geomtri: exchanging two flags}
\[
    \begin{array}{rl}
      \geomtri_1(\F_1,\F_3,\F_2)&=-\geomtri_1(\F_1,\F_2,\F_3),\\  
      \geomtri_2(\F_1,\F_3,\F_2)&=-\geomtri_3(\F_1,\F_2,\F_3)
    \end{array}
    \; ;\]

  \item
    \label{i- sum geomtri is 0}
 We have $\geomtri_1(\Tau)+\geomtri_2(\Tau)+\geomtri_3(\Tau)=0$;

  \item  
    \label{i- geomtri is ultram} 
    For all $\indTriRap\in\ZZ/3\ZZ$,
    if $\geomtri_\indTriRap(\Tau)>0$,
    then we have $\geomtri_{\indTriRap-1}(\Tau)=0$ and 
    $\geomtri_{\indTriRap+1}(\Tau)=-\geomtri_\indTriRap(\Tau)<0$.

  \end{myenum}
\end{proposition}

In order to prove this proposition, in particular,
the invariance of the triple ratio by cyclic
permutation of the flags,
we first introduce 
the natural dual invariants given by the cross ratios 
of the natural induced quadruple of points on the line $\D_1$ 
(that is, exchanging the role of points and lines):
\[    \begin{array}{rl}
\geomtri^*_1(\F_1,\F_2, \F_3) 
&=  \geombir(\p_1,\D_2 \cap \D_1, \D_{23} \cap \D_1, \D_3 \cap \D_1)\\
\geomtri^*_2(\F_1,\F_2,\F_3)
&=  \geombir(\p_1, \D_3\cap\D_1, \D_2 \cap \D_1, \D_{23} \cap \D_1)\\
\geomtri^*_3(\F_1,\F_2,\F_3)
&= \geombir(\p_1, \D_{23} \cap \D_1, \D_3\cap\D_1, \D_2 \cap \D_1)\;.
    \end{array}
\]

The following property is straigthforward.
\begin{equation}
  \label{eq- geomtriet: exchanging two flags}
  \begin{array}{rl}
    \geomtri^*_1(\F_1,\F_3,\F_2)&=-\geomtri^*_1(\F_1,\F_2,\F_3)\\  
    \geomtri^*_2(\F_1,\F_3,\F_2)&=-\geomtri^*_3(\F_1,\F_2,\F_3)\;.
  \end{array}
\end{equation}

We will need the following property showing
that the invariants behave nicely under duality.

\begin{lemma}
\label{lemm- geomtriet and geomtri}
For $\indTriRap=1,2,3$, we have 
$\geomtri^*_\indTriRap(\F_1,\F_2, \F_3)
=  \geomtri_\indTriRap(\F_3,\F_2,\F_1)$.
\end{lemma}

\begin{proof}[Proof of Lemma \ref{lemm- geomtriet and geomtri}]
By invariance under perspectivities and double transpositions,
we have
\begin{align*}
\geomtri^*_1(\F_1,\F_2,\F_3) 
&=  \geombir(\p_1, \D_2\cap\D_1, \D_{23}\cap\D_1, \D_3\cap\D_1)\\
&=  \geombir(\p_1\p_3, \p_{12}\p_3, \D_{23}, \D_3)\\
&=  \geombir(\D_3, \p_2\p_3, \p_{12}\p_3, \p_1\p_3)\\
&= \geomtri_1(\F_3,\F_2, \F_1)\;.
\end{align*}
The proof of $\geomtri^*_\indTriRap(\F_1,\F_2, \F_3)=
\geomtri_\indTriRap(\F_3,\F_2,\F_1)$ for $\indTriRap=2,3$ is similar.
\end{proof}

We now turn to the proof of Proposition \ref{prop- props of geomtri}. 

\begin{proof}[Proof of Proposition \ref{prop- props of geomtri}]
Assertions \ref{i- sum geomtri is 0} 
and \ref{i- geomtri is ultram}
follow immediately 
from the properties of the cross ratio $\geombir$ 
under cyclic permutation of the three last points
(see \eqref{eq- sum of 3-cyclic perm of geom bir is 0} 
and
 \eqref{eq- geombir is ultram}).

Assertion \ref{i- geomtri: exchanging two flags} 
follows immediately from the definition and from the symmetries of
the cross ratio.

We finally 
prove \ref{i- geomtri: invariant by cyclic  perm}
of Proposition \ref{prop- props of geomtri}.
Using \ref{i- geomtri: exchanging two flags}, 
Lemma \ref{lemm- geomtriet and geomtri}
and
\eqref{eq- geomtriet: exchanging two flags},
we have
  \begin{align*}
    \geomtri_1(\F_2,\F_3,\F_1)
    &=-\geomtri_1(\F_2,\F_1,\F_3)\\ 
    &=-\geomtri^*_1(\F_3,\F_1, \F_2)\\
    &=\geomtri^*_1(\F_3,\F_2, \F_1)
    =\geomtri_1(\F_1,\F_2,\F_3) ,
  \end{align*}
  \begin{align*}
    \geomtri_2(\F_2,\F_3,\F_1)
    &=-\geomtri_3(\F_2,\F_1,\F_3)\\ 
    &=-\geomtri^*_3(\F_3,\F_1, \F_2)\\
    &=\geomtri^*_2(\F_3,\F_2, \F_1)
    =\geomtri_2(\F_1,\F_2,\F_3)\;.
  \end{align*}
The case where $\indTriRap=3$ is similar to the case $\indTriRap=2$.
\end{proof}

\subsection{Geometric triple ratio from  algebraic triple ratio}
When $\PP$ is the projective plane  
on some field $\KK$  endowed with some ultrametric absolute value, 
 and $\geombir = \log\abs{\Bir}$ where $\Bir$ 
is the usual $\KK$-valued cross ratio on $\PP$,
the three geometric triple ratios $\geomtri_\indTriRap(\Tau)$, $\indTriRap=1,2,3$
of $\Tau$
 are obtained from 
the single algebraic triple ratio $\Z=\Tri(\Tau)$ of $\Tau$
 by  the relations 
\begin{equation}
  \begin{array}{rcl}
    \geomtri_1(\Tau)&=&\log\abs{\Z}\\
    \geomtri_2(\Tau) &=&\log\abs{\frac{1}{1+\Z}}=-\log\abs{1+\Z}\\
    \geomtri_3(\Tau) &=& \log\abs{1+\Z^{-1}},
  \end{array}
\end{equation}
which are easily derived from the expression of 
algebraic triple ratio as a  cross ratio \eqref{eq- triple ratio as a cross ratio}
and
from the symmetry properties of the algebraic cross ratio  \eqref{eq- sym bir 1342}.

\begin{remark}
Note that the geometric invariants do not determine 
the triple of flags up to automorphisms of $\PP$
(unlike the usual (algebraic) triple ratio): 
for example in the algebraic case  $\PP=\PP(\KK^3)$,  
take $\Tau$ with triple ratio $\Z\in\KK$ with $\abs{\Z}>1$ 
and $\Tau'$ with triple ratio $\Z'=\Z\aK$ where $\aK\in\KK$ with $\abs{\aK}=1$ and
$\aK\neq 1$. Then $\Tau$ and $\Tau'$ are not in the same
$\PGL(\KK^3)$-orbit, but have the same three geometric
invariants, as
$\geomtri_1(\Tau) = \log\abs{\Z}=\geomtri_1(\Tau')$, 
$\geomtri_2(\Tau) =-\log\abs{\Z}=\geomtri_2(\Tau')$,
$\geomtri_3(\Tau)= 0=\geomtri_3(\Tau')$.
\end{remark}

\section{Proof of the main result}
In this section we
prove Theorems \ref{theo- tripod} and \ref{theo- triangle}.
Let $\Es$ be a Euclidean
building of type $A_2$ and 
$\Tau=(\F_1,\F_2,\F_3)$  be a \genericTF{} triple of flags 
in the projective plane $\PP$ at infinity of $\Es$.
We denote by 
$\gZ_\indTriRap=\geomtri_\indTriRap(\F_1,\F_2,\F_3)$, 
$\indTriRap=1,2,3$,
its geometric triple ratio,
and by $\Aij=\App(\F_\ii,\F_\jj)$, $\Ap=\App(\p_1,\p_2,\p_3)$ and
$\AD=\App(\D_1,\D_2,\D_3)$ the five associated flats.

We first define the six associated points in $\Es$.

 \subsection{Associated  points in the building}
For $\{\ii,\jj,\kk\}=\{1,2,3\}$,
denote by $\yTauE_\kk$  the center in $\Es$ of the projective frame
$(\p_1,\p_2,\p_3, \p_{\ii\jj})$, where $\p_{\ii\jj}=\D_\ii \cap \D_\jj$,
and by $\yTauE_\kk^*$   the center of the projective frame
$(\D_1,\D_2,\D_3, \D_{\ii\jj})$, where $\D_{\ii\jj}=\p_\ii\p_\jj$,
as defined in 
Proposition \ref{prop- center of a projective frame}.
In particular
the point $\yTauE_\kk$ is the orthogonal projection of $\p_{\ii\jj}$ on $\Ap$,
the point $\yTauE_\kk^*$ is the orthogonal projection of $\D_{\ii\jj}$ on  $\AD$,
the point $\yTauE_\kk$ is the  orthogonal projection of $\p_\kk$ 
on $\Aij=\App(\p_\ii,\p_\jj,\p_{\ii\jj})$, 
and
the point $\yTauE_\kk^*$ is the orthogonal projection of $\D_\kk$ 
on $\Aij=\App(\D_\ii,\D_\jj,\D_{\ii\jj})$.

\subsection{In the flat \texorpdfstring{$\Aij$}{A\_ij}}
We now link 
the respective position of the points  $\yTauE_\kk$ and $\yTauE_\kk^*$ in the flat
$\Aij$ to the geometric triple ratio of $\Tau$. 
Suppose that the indices $\ii,\jj,\kk$ respect the cyclic order, 
i.e. that $(\ii,\jj,\kk)=(1 2 3)$ as cyclic permutations.
We identify $\Aij$ with the model flat $\Aa$ 
by a marked flat $\f_{\ii\jj}:\Aa \to \Aij$ sending $\bord\Cc$ to
$\F_\jj$.
For $\xE,\yE$ in $\Aij\simeq\Aa$, 
we define then $\ovra{\xE\yE}=\yE-\xE=\Bus{\F_\jj}{\xE}{\yE}$.
Recall that 
$(\eRRNo_1,\eRRNo_2,\eRRNo_3)$ denotes the canonical basis of $\RR^3$.
In particular, the  directions of $\p_\ii$, $\p_{\ii\jj}$ and $\p_\jj$  
are respectively identified with the directions of 
$\classA{\eRRNo_1}$, $\classA{\eRRNo_2}$, and  $\classA{\eRRNo_3}$ in $\Aa$.

\begin{proposition}
\label{prop- in Aij - y*y}
The following holds.
\begin{myenum}
\item 
\label{i- y*y in simple roots coords} 
In simple roots coordinates, we have
$\ovra{\yTauE_\kk^* \yTauE_\kk}= (\gZ_2, \gZ_3)$;

\item 
\label{i- yy* sing of type point} 
For $\indTriRap=1,2,3$, if $\gZ_\indTriRap>0$
then 
$\ovra{\yTauE_\kk\yTauE_\kk^*}=\gZ_\indTriRap
\classA{\eRRNo_\indTriRap}$. 
In particular
$\yTauE_\kk^*$ is on one of the three singular rays of type
  $\stypep$ issued from $\yTauE_\kk$ 
(i.e the rays to $\p_\ii$,  $\p_\jj$ and $\p_{\ii\jj}$).
\end{myenum}
\end{proposition}

\begin{proof}
As $\yTauE_\kk$ and $\yTauE_\kk^*$ are the respective orthogonal projections on the flat
$\App_{\ii\jj}$ of $\p_\kk$ and $\D_\kk$, 
by Proposition \ref{prop- 1 flat p D} and  cyclic invariance of the geometric triple ratio, we have
\begin{align*}
&
\rac1( \ovra{\yTauE_\kk^* \yTauE_\kk})
=  \geombir(\D_\ii, \p_\ii\p_\kk,  \p_\ii\p_\jj , \p_\ii \p_{\jj\kk} )
= \geomtri_2(\F_\ii, \F_\jj, \F_\kk)
= \gZ_2\\
\mbox{and }
&
\rac2( \ovra{\yTauE_\kk^* \yTauE_\kk})
=  \geombir(\D_\jj, \p_\jj \p_{\kk\ii},  \p_\jj\p_\ii , \p_\jj\p_\kk )
= \geomtri_3(\F_\jj, \F_\kk, \F_\ii)=\gZ_3
\;.
\end{align*}

Assertion \ref{i- yy* sing of type point} 
follows, since we have then 
$\gZ_{\indTriRap-1}=0$ and $\gZ_{\indTriRap+1}=-\gZ_\indTriRap$ 
by ultrametricity of the geometric triple ratio
(Proposition \ref{prop- props of geomtri}\ref{i- geomtri is
  ultram}).
\end{proof}

We now describe the intersections of $\Aij$ with the four other flats
(see Figures \ref{fig-  type tripod} and 
\ref{fig-  type triangle in positive case} in the introduction). 
These intersections happen to be  sectors in $\Aa$ bounded by two
singular rays of same type, equivalently the union of two adjacent
Weyl chambers.

\begin{proposition}
\label{prop- in Aij - intersections}
Let $\xE\in\Aij$. Then

\begin{myenum}

\item
\label{i- Aij cap Ap}
The intersection $\Aij\cap\Ap$ is 
the sector 
at $\yTauE_\kk$  
bounded by the rays to $\p_\ii$ and $\p_\jj$.
That is
\[\xE\in  \Ap
\mbox{  if and only if }
\left\{
\begin{array}{lll}
\rac1(\xE) & \geq & \rac1(\yTauE_\kk) \\
\rac2(\xE) & \leq & \rac2(\yTauE_\kk)
\end{array}
\right.
\;.\]

\item
\label{i- Aij cap AD}
The intersection $\Aij\cap\AD$ is 
the sector at $\yTauE_\kk^*$  
bounded by the rays to $\D_\ii$ and $\D_\jj$.
That is, 
\[\xE\in  \AD
\mbox{  if and only if }
\left\{
\begin{array}{lll}
\rac1(\xE) & \leq & \rac1(\yTauE_\kk^*) \\
\rac2(\xE) & \geq & \rac2(\yTauE_\kk^*)
\end{array}
\right.
\;.\]

\item
\label{i- Aij cap Ajk}
The intersection $\Aij\cap\Ajk$ is the intersection of 
the sector at $\yTauE_\kk$ bounded by the rays to $\p_\jj$ and
$\D_\ii\cap\D_\jj$,
and
the sector at $\yTauE_\kk^*$ bounded by the rays to $\D_\jj$ and
$\p_\ii\p_\jj$. 
That is,
\[\xE\in  \Ajk
\mbox{  if and only if }
\left\{
\begin{array}{lll}
\rac1(\xE) & \geq &  \rac1(\yTauE_\kk^*) \\
\rac2(\xE) & \geq & \rac2(\yTauE_\kk) \\
\rac3(\xE) & \leq & \min(\rac3(\yTauE_\kk), \rac3(\yTauE_\kk^*))
\end{array}
\right.
\;.\]

\item
\label{i- Aij cap Aki}
The intersection $\Aij\cap\Aki$ is the intersection of 
the sector at $\yTauE_\kk$ bounded by the rays to $\p_\ii$ and
$\D_\ii\cap\D_\jj$,
and
the sector at $\yTauE_\kk^*$ bounded by the rays to $\D_\ii$ and
$\p_\ii\p_\jj$.
That is,
\[\xE\in  \Aki
\mbox{  if and only if }
\left\{
\begin{array}{lll}
\rac1(\xE) & \leq &  \rac1(\yTauE_\kk) \\
\rac2(\xE) & \leq & \rac2(\yTauE_\kk^*) \\
\rac3(\xE) & \geq & \max(\rac3(\yTauE_\kk), \rac3(\yTauE_\kk^*))
\end{array}
\right.
\;.\]
\end{myenum}
\end{proposition}

\begin{proof}
Since $\yTauE_\kk$ is the center of the projective frame
$(\p_\ii,\p_\jj, \p_{\ii\jj},\p_\kk)$, 
assertion \ref{i- Aij cap Ap} comes from
Proposition \ref{prop- Ai cap Aj and partition Ai}, as
$\App_{\ii\jj}=\App(\p_\ii,\p_\jj, \p_{\ii\jj})$
and  $\Ap=\App(\p_\ii,\p_\jj,\p_\kk)$.
Assertion \ref{i- Aij cap AD} is similar. 
Assertion \ref{i- Aij cap Ajk}:
A point $\xE\in\Aij$ lies in $\Ajk$ 
if and only if, 
in the spherical building of directions at $\TangS_\xE\Es$, 
the direction $\TangS_\xE\D_\jj$ is  opposite to $\TangS_\xE\p_\kk$  
and  $\TangS_\xE\p_\jj$ is  opposite to $\TangS_\xE\D_\kk$.
Moreover, 
$\TangS_\xE\D_\jj$ is  opposite to $\TangS_\xE\p_\kk$
if and only if
$\xE\in\App(\p_\kk,\p_\jj,\p_{\ii\jj})$.
As $\yTauE_\kk$ is the center of the projective frame 
$(\p_\ii,\p_\jj,\p_{\ii\jj}, \p_\kk)$ and 
$\App_{\ii\jj}=\App(\p_\ii,\p_\jj,\p_{\ii\jj})$,
the set of such $\xE$ is
the sector at $\yTauE_\kk$ bounded by 
the rays to $\p_\jj$ and $\D_\ii\cap\D_\jj$
(by Proposition \ref{prop- Ai cap Aj and partition Ai}).
This is the subset of $\xE\in\Aij$ satisfying:
$\rac2(\xE)\geq \rac2(\yTauE_\kk)$ and
$\rac3(\xE)\leq\rac3(\yTauE_\kk)$.
Similarly,
as $\yTauE^*_\kk$ is the center of the projective frame 
$(\D_\ii,\D_\jj,\D_{\ii\jj}, \D_\kk)$ and 
$\App_{\ii\jj}=\App(\D_\ii,\D_\jj,\D_{\ii\jj})$,
the direction $\TangS_\xE\p_\jj$ is  opposite to $\TangS_\xE\D_\kk$
if and only if
$\xE$ is in the sector 
at $\yTauE_\kk^*$ 
bounded by the rays to $\D_\jj$ and $\D_{\ii\jj}=\p_\ii\p_\jj$.
That is, if and only if
$\rac1(\xE)\geq \rac1(\yTauE_\kk^*)$ and
$\rac3(\xE)\leq\rac3(\yTauE_\kk^*)$, and we are done.
Assertion \ref{i- Aij cap Aki} is similar.
\end{proof}

In particular, as $\yTauE_\kk^*$ is on one of the three singular rays of type
  $\stypep$ issued from $\yTauE_\kk$ by Propositions \ref{prop- in Aij - y*y}, 
from Proposition \ref{prop-  in Aij - intersections} we easily 
get the following result.

\begin{corollary}
\label{coro- partition Aij}
The intersections with $\Aij$ of $\Ajk$,$\Aki$, $\Ap$ and $\AD$ 
form a partition  of $\Aij$.  
\qed\end{corollary}

\subsection{In the flat \texorpdfstring{$\Ap$}{A\_p}}
We now consider the flat $\Ap=\App(\p_1,\p_2,\p_3)$.
The following proposition describes the respective positions in $\Ap$
of the points $\yTauE_1$, $\yTauE_2$, $\yTauE_3$.
We identify $\Ap$ with $\Aa$ by a marked flat $\f_\p:\Aa \to
\Ap$ sending $\bord\Cc$ to $(\p_1,\p_1\p_2)$
(hence direction $\classA{\eRRNo_\ii}$ to $\p_\ii$ for $\ii=1,2,3$). 
Recall that we then have
$\ovra{\xE\xE'}=\xE'-\xE=\Bus{(\p_1,\p_1\p_2)} {\xE}{\xE'}$  for
$\xE,\xE'\in\Ap$.

\begin{proposition}
\label{prop- in Ap}
In the flat $\Ap$ we have:

 \begin{myenum}
\item
 \label{i- y2 y3}
In simple roots coordinates, we have $\ovra{\yTauE_2\yTauE_3}=(\gZ_1,0)$.

 \item 
\label{i- y_i y_i+1}
If $\gZ_1\geq 0$,
the point  $\yTauE_{\ii+1}$ is in the ray $[\yTauE_\ii,\p_{\ii+2})$ 
(for all $\ii$),
and if $\gZ_1\leq 0$,
the point $\yTauE_\ii$ is in the ray $[\yTauE_{\ii+1},\p_{\ii+2})$ for
all $\ii$.

\end{myenum}

In particular the triangle $\Delta\subset \Ap$ 
with vertices $\yTauE_1$, $\yTauE_2$,  $\yTauE_3$ 
is {\em  singular}, 
i.e. the sides have singular type in $\Cb$.
\end{proposition}

\begin{proof}
Recall that the point $\yTauE_\kk$ is the orthogonal projection 
on the flat $\Ap$ 
of the singular boundary point $\p_{\ii\jj}=\D_\ii \cap \D_\jj$.
Then, by Proposition \ref{prop- proj of centers to transverse spaces}
the points  $\yTauE_2$ and $\yTauE_3$ have the same projection 
in the transverse tree $\Es_{\p_1}$,
that is the center of the ideal triple
$ (\p_1\p_{13}, \p_1\p_2, \p_1\p_3)
= (\D_1, \p_1\p_2, \p_1\p_3)
= (\p_1\p_{23}, \p_1\p_2, \p_1\p_3)
$,
proving that $\rac2(\ovra{\yTauE_2\yTauE_3})=0$.
Furthermore, 
by Proposition \ref{prop- 1 plat, 2 points} 
we have
\begin{align*}
  \rac2(\ovra{\yTauE_3 \yTauE_1})
  &= \geombir(\p_1\p_2, \p_1\p_{12},  \p_1\p_3 , \p_1\p_{23} )\\
  &= \geombir(\p_1\p_2, \D_1,  \p_1\p_3 , \p_1\p_{23} )\\
  &= \geombir(\D_1, \p_1\p_2, \p_1\p_{23}, \p_1\p_3)\\
  &=\gZ_1
\end{align*}
proving that $\rac2(\ovra{\yTauE_3\yTauE_1})=\gZ_1$.
Applying this to the permuted triple $(\F_3,\F_1,\F_2)$, we obtain
$\rac1(\ovra{\yTauE_2\yTauE_3}) = \gZ_1$ (by invariance of the
geometric triple
ratio $\gZ_1$ by cyclic permutation).
Assertion \ref{i- y_i y_i+1}
follows from \ref{i- y_i y_i+1}, applying cyclic permutations.
\end{proof}

We now describe the intersections of $\Ap$ with the other flats, 
see Figure \ref{fig-  in Ap}.

\begin{figure}[h]
\centering
\hfill
  \begin{minipage}[b]{0.4\linewidth}
\centering
\includegraphics[scale=0.8]{figures/Es_Ap_triangle_pos_yTauE_avec_gZ.fig}

In the case $\gZ_1\geq 0$.
  \end{minipage}
\hfill
\begin{minipage}[b]{0.4\linewidth}
\centering
\includegraphics[scale=0.8]{figures/Es_Ap_triangle_neg_yTauE_avec_gZ.fig}

In the case $\gZ_1\leq 0$.
  \end{minipage}
\hfill

\caption{In the flat $\Ap$.}
\label{fig-  in Ap}
\end{figure}
\begin{proposition}
\label{prop- in Ap - intersections}
Let $\SW_\ii=\Ap \cap  \App_{\ii,\ii+1}$ 
and let $\Delta$ be the 
triangle with vertices $\yTauE_1,\yTauE_2,\yTauE_3$.
Then

\begin{myenum}
\item 
\label{i- Ap cap Aij}
$\SW_\ii$ is the sector of $\Ap$ bounded by the rays 
from $\yTauE_{\ii+2}$ to $\p_\ii$ and $\p_{\ii+1}$.

\item 
\label{i- partition Ap}
$\SW_1,\SW_2,\SW_3$ and $\Delta$ form a partition of $\Ap$.

\end{myenum}
\end{proposition}

\begin{proof}
Assertion \ref{i- Ap cap Aij} follows from
point \ref{i- Aij cap Ap} 
of Proposition \ref{prop- in Aij - intersections}.
In the case where $\gZ_1\geq 0$, 
assertion \ref{i- partition Ap} then comes 
from the fact that 
for all $\ii$, $\yTauE_{\ii+1}$ is in the ray $[\yTauE_\ii,\p_{\ii+2})$ 
(Proposition \ref{prop- in Ap}).
The case where $\gZ_1\leq 0$ is similar.
\end{proof}

\subsection{In the flat \texorpdfstring{$\AD$}{A\_D}}
We now state the similar properties in the dual  flat $\AD=\App(\D_1,\D_2,\D_3)$,
which have same proofs, exchanging the role of points and lines.

\begin{proposition}
\label{prop- in AD}
In the flat $\AD$  identified with $\Aa$ by a marked flat 
sending $\bord\Cc$ to $(\D_1\cap\D_2,\D_1)$, we have:
\begin{myenum}
\item   $\ovra{\yTauE^*_2\yTauE^*_3}=(0,-\gZ_1)$
in simple roots coordinates.
In particular $\yTauE^*_2$ and $\yTauE^*_3$ 
are on a common singular geodesic to $\D_1$.
\item The points  $\yTauE^*_1$, $\yTauE^*_2$,  $\yTauE^*_3$ 
form a singular triangle $\Delta^*$ in $\AD$.

\item 
  \label{i- AD cap Aij}
  For all $\ii\in\ZZ/3\ZZ$, $\SW^*_\ii=\AD \cap \App_{\ii,\ii+1}$ is
  the sector of $\AD$ bounded by the rays from $\yTauE^*_{\ii+2}$
  to $\D_\ii$ and $\D_{\ii+1}$.

\item 
  \label{i- partition AD}
  $\SW^*_1,\SW^*_2,\SW^*_3$ and $\Delta^*$ form a partition of $\AD$.

\end{myenum}
\qed
\end{proposition}

\begin{figure}[h]
\centering
\hfill
  \begin{minipage}[b]{0.4\linewidth}
\centering
\includegraphics[scale=0.8]{figures/Es_AD_triangle_pos_yTauE.fig}
    
In the case $\gZ_1\geq 0$.
  \end{minipage}
\hfill
  \begin{minipage}[b]{0.4\linewidth}
\centering
\includegraphics[scale=0.8]{figures/Es_AD_triangle_neg_yTauE.fig}
    
In the case $\gZ_1\leq 0$.
  \end{minipage}
\hfill

\caption{In the flat $\AD$.}
\label{fig-  in AD}
\end{figure}

\subsection{The classification}
We now combine the previous results  
to establish the classification in two geometric types, 
finishing to prove Theorems \ref{theo- tripod} and \ref{theo- triangle}.

\begin{proof}[Proof of Theorem \ref{theo- tripod}]
Let $\vstE=\yTauE_3$ and $\vstE^*=\yTauE_3^*$.
We identify the flat $\App_{12}$
with the model flat $\Aa$ 
by a marked flat sending $\bord\Cc$ to $\F_2$, and $0$ to 
$\yTauE_3^*$.
By Proposition \ref{prop- in Aij - intersections}
applied to the flat $\App_{12}$, we have
$\rac1(\yTauE_3)=\gZ_2$,
$\rac2(\yTauE_3)=\gZ_3$, and
$\rac3(\yTauE_3)=\gZ_1$.
By Proposition \ref{prop- in Aij - intersections}
applied to the flat $\App_{12}$,
the intersection $I=\App_{12}\cap \App_{23} \cap \App_{31}$ 
is the  subset of $\yE\in \App_{12}$ such that
\[
\left\{
\begin{array}{lll}
  0 &\leq \rac1(\yE) &\leq \rac1(\yTauE_3)=\gZ_2\\  
  0 &\geq \rac2(\yE) &\geq  \rac2(\yTauE_3)=\gZ_3\\ 
  \max(\rac3(\yTauE_3), 0)& 
  \leq \rac3(\yE)  
  &\leq  \min(\rac3(\yTauE_3), 0) 
  \; .
\end{array}
\right.
\]
In particular, 
if $I$ is not empty, then
$\gZ_1=\rac3(\yTauE_3)=0$.

Suppose from now on that $\gZ_1=0$. 
Then  $\gZ_2\geq 0$ and $\gZ_3 =-\gZ_2$  by 
the ultrametricity of the geometric triple ratio
(Proposition \ref{prop- props of geomtri}\ref{i- geomtri is ultram}).
By the description above, $I$ 
is then the subset of the line $\rac3=0$ (which contains
$\yTauE_3^*=0$ and $\yTauE_3$) consisting of the $\yE$ such that
$ 0 \leq \rac1(\yE) \leq \rac1(\yTauE_3)$
(since $\rac2(\yE)=-\rac1(\yE)$ when $\rac3(\yE)=0$).
Hence $I$ is  not empty and 
is the segment from $0=\yTauE_3^*$ to $\yTauE_3$ 
i.e. $[\vstE,\vstE^*]$.
Furthermore, 
as $\gZ_1=0$,  Proposition \ref{prop- in Ap}  implies that 
$\yTauE_1=\yTauE_2=\yTauE_3$. 
Similarly, 
we have  $\yTauE_1^*=\yTauE_2^*=\yTauE_3^*$ 
by Proposition \ref{prop- in AD}.
Suppose now $\vstE\neq \vstE^*$.
Since the segment $[\vstE,\vstE^*]$ lies in the ray
$[\vstE,\p_{\ii\jj})$,
and $\vstE=\yTauE_\kk$ is the orthogonal projection of
$\p_{\ii\jj}$ on $\Ap$, 
we have $\sphericalangle_\vstE(\vstE^*,\D)=\pi$ 
for all lines $\D$ in $\bordinf\Ap$ 
(Proposition \ref{prop- center of a projective frame}).
Therefore  we have
$\sphericalangle_\vstE(\vstE^*,\yE)\geq \frac{2\pi}{3}$
for all $\yE\neq\vstE$ in $\Ap$.
Similarly,  we have that $\sphericalangle_{\vstE^*}(\vstE,\yE)\geq
\frac{2\pi}{3}$ for all $\yE\neq\vstE$ in $\Ap$.
Hence  $[\vstE,\vstE^*]$ is the unique segment of minimal length joining
$\Ap$ to $\AD$.
Assertion \ref{i- coords xx^*} 
follows from Proposition \ref{prop- in  Aij - y*y}.
\end{proof}

\begin{proof}[Proof of Theorem \ref{theo- triangle}]
If $\gZ_2 > 0$, then  $\gZ_1 = 0$
by the ultrametricity of the geometric triple ratio
(Proposition \ref{prop- props of geomtri}\ref{i- geomtri is ultram}), 
and  $\Ap \cap \AD$ is empty by 
Theorem \ref{theo- tripod}.
Suppose now that $\gZ_2\leq 0$. 
Since the case $\gZ_1 \leq 0 $ reduces to the case $\gZ_1 \geq 0 $ by 
exchanging $\F_2$ and $\F_3$,
it is enough to handle the case $\gZ_1 \geq 0$. 
Then $\gZ_3=0$ and $\gZ_2=-\gZ_1$. 
Let $\vstE_\ii=\yTauE_{\ii+2}$ for $\ii\in\ZZ/3\ZZ$.
In $\Aij$ identified with $\Aa$ in such a way that $\yTauE_\kk^*=0$,
by Proposition \ref{prop- in Aij - y*y} we have 
$\rac1(\yTauE_\kk)=\gZ_2=-\gZ_1 \leq 0$,
 $\rac2(\yTauE_\kk)=\gZ_3 = 0$,
hence $\rac3(\yTauE_\kk)=\gZ_1 \geq 0$.
By Proposition \ref{prop- in Aij -  intersections}\ref{i- Aij cap Aki},
$\Aij\cap\Aik$ is the set of $\xE\in \Aij\simeq\Aa$
such that
$\rac1(\xE)\leq   \rac1(\yTauE_\kk)$,
$\rac2(\xE)  \leq  0=\rac2(\yTauE_\kk)$
and $\rac3(\xE)  \geq  \max(\rac3(\yTauE_\kk), 0)=\rac3(\yTauE_\kk)$.
This is the Weyl chamber $\yTauE_\kk-\Cb$, 
i.e. the Weyl chamber from  $\yTauE_\kk=\vstE_\ii$ to $\F_\ii$.
Similarly, $\Aij\cap\Ajk$ 
is the Weyl chamber from  $\yTauE_\kk^*$ to $\F_\jj$.
Applying a cyclic permutation $(\ii \jj \kk)$,
i.e. working in the flat $\Ajk$, 
we also similarly  get that $\Aij\cap\Ajk$ 
is the Weyl chamber from  $\yTauE_\ii$ to $\F_\jj$.
Therefore $\yTauE_\kk^*=\yTauE_\ii$.

By Proposition \ref{prop- in Aij - intersections}
$\Ap \cap \AD \cap \Aij$ 
is the intersection 
of the sector at $\yTauE_\kk^*$  
bounded by the rays to $\D_\ii$ and $\D_\jj$,
with the sector  at $\yTauE_\kk$  
bounded by the rays to $\p_\ii$ and $\p_\jj$.
As the point $\yTauE_\kk$ is on the ray from $\yTauE_\kk$ to $\D_\ii$, 
this is equal to the segment $[\yTauE_\kk,\yTauE_\kk^*]$. 
In particular $\Ap \cap \AD $ contains $\yTauE_\kk$.
Then
$\Ap \cap \AD$ contains $\yTauE_1$, $\yTauE_2$ and $\yTauE_3$, hence
the triangle $\Delta$ with vertices $\yTauE_1$, $\yTauE_2$ and $\yTauE_3$, 
and  
since 
$\Ap \cap \AD \cap \Aij 
= [\yTauE_\kk,\yTauE_\ii]
\subset \Delta$, 
Proposition \ref{prop- in  Ap - intersections}\ref{i- partition Ap} 
provides the reverse inclusion. 
Assertion \ref{i- xi xj} comes from Proposition 
\ref{prop- in Aij -  y*y}.

We finally prove \ref{i- Delta opp Fi en xi}. 
Let $(\ii,\jj,\kk)=(1 2 3)$.
Looking in the flat $\Ap$, we see that 
the singular triangle $\Delta$ is contained in the Weyl chamber of $\Es$
with tip $\vstE_\ii$  and that
at $\vstE_\ii$, we have
$\TangS_{\vstE_\ii} \vstE_\jj = \TangS_{\vstE_\ii} \p_\jj$.
Looking in the flat $\AD$ we get
$\TangS_{\vstE_\ii} \vstE_\kk = \TangS_{\vstE_\ii} \D_\kk$. 
Hence $\TangS_{\vstE_\ii} \Delta=(\TangS_{\vstE_\ii}\p_\jj,\TangS_{\vstE_\ii} \D_\kk)$.
Since
$\vstE_\ii$ belongs to the flats $\App(\F_\ii,\F_\jj)$ and
$\App(\F_\ii,\F_\kk)$, 
we have that
$\TangS_{\vstE_\ii} \p_\jj$ is opposite to $\TangS_{\vstE_\ii} \D_\ii$
and that
$\TangS_{\vstE_\ii} \D_\kk$ is opposite to $\TangS_{\vstE_\ii} \p_\ii$.
Therefore the Weyl chambers
$\TangS_{\vstE_\ii} \Delta$ and $\TangS_{\vstE_\ii} \F_\ii$ are opposite.
It implies that $\Delta$ and the Weyl chamber from $\vstE_\ii$ to
$\F_\ii$ are contained in a common flat of $\Es$ by basic properties of
real Euclidean buildings (see property (CO) of
\cite{ParImm}).
\end{proof}

In the algebraic case the following remark 
provides an alternative proof of some of the assertions 
of  Theorem  \ref{theo- triangle}.
\begin{remark}
Let $\pt_\ii$ in $\V=\KK^3$ be a vector  representing $\p_\ii$ and
$\Dt_\ii$ in $\V^*$ be a linear form  representing $\D_\ii$.
Let $\eVb=(\eV_1,\eV_2,\eV_3)$ 
be the basis of $\V$ dual to the basis
$(\Dt_1,\Dt_2,\Dt_3)$ of $\V^*$. 
Then in the projective plane $\classPV{\eV_\ii}=\D_\jj\cap \D_\kk$.
We may suppose that $\pt_1=(0,1,1)$, $\pt_2=(\Z,0,1)$, $\pt_3=(1,1,0)$
in the basis $\eVb$, with $\Z=\Tri(\F_1,\F_2,\F_3)$.
Then the element $\g\in\GL(V)$ with matrix in the basis $\eVb$ 
\[
\begin{pmatrix}
1   & 1 & 0\\
0   & 1 & 1 \\
1/\Z & 0 & 1
\end{pmatrix}\]
sends $\classPV{\eV_\ii}$ to $\p_{\ii+1}$, 
hence $\AD$ to $\Ap$.
If $\abs{1+\Z}\geq 1$ and
$\gZ=\log\abs{\Z}\geq 0$, 
then the  fixed point set  of $\g$ in $\AD$ 
is the image by the marked flat $\f_\eVb$
of the singular triangle 
$\{\vA\in\Cb \ | \ \vA_1-\vA_3 \leq \log\abs{\Z}\}$ (that is, $\Delta$).
\end{remark}

\subsection{Complements}
We add here for future use  a simple description
 of the vertices $\vstE_\ii,\vstE_\jj,\vstE_\kk$ of the singular triangle  $\Delta$
in Theorem \ref{theo- triangle}
by the projections on transverse trees at infinity.

\begin{lemma}
\label{lemm- proj de x_i sur les arbres transverses}
Under the hypotheses and notations of Theorem \ref{theo- triangle},
we have the following properties.
\begin{myenum}
\item 
\label{i- proj de xEi sur tree at p_i}
The projection $\projE_{\p_\ii}(\vstE_\ii)$ of $\vstE_\ii$ on the tree $\Es_{\p_\ii}$  
is  the center of the ideal tripod $(\D_\ii, \p_\ii\p_\jj, \p_\ii\p_\kk)$.

\item 
\label{i- proj de xEi sur tree at D_i}
The projection $\projE_{\D_\ii}(\vstE_\ii)$ of $\vstE_\ii$  on the tree
$\Es_{\D_\ii}$
is  the center  of the ideal tripod $(\p_\ii, \D_\ii\cap\D_\jj, \D_\ii\cap\D_\kk)$.

\item
\label{i- proj de xEj sur tree at p_i}
 The projection $\projE_{\p_\ii}(\vstE_\jj)$ 
is  the center of the ideal tripod 
$(\D_\ii, \p_\ii\p_\jj, \p_\ii\p_{\jj\kk})$.

\item
\label{i- proj de xEj sur tree at D_i}
 The projection  $\projE_{\D_\ii}(\vstE_\jj)$ 
is  the center of the ideal tripod 
$(\p_\ii, \D_\ii\cap\D_\jj, \D_\ii\cap \D_{\jj\kk})$.
\end{myenum}

\end{lemma}

\begin{proof}
As the point $\vstE_\ii$ belongs to the three flats 
$\App(\F_\kk,\F_\ii)$ and $\App(\F_\jj,\F_\ii)$ and $\App(\p_\ii,\p_j,\p_k)$,
its projection in the tree $\Es_{\p_\ii}$
belongs to the projection of $\App(\F_\jj,\F_\ii)$, 
which is the line from $\D_\ii$ to $\p_\ii\p_\jj$,
to the projection of $\App(\F_\kk,\F_\ii)$, 
which is the line from $\D_\ii$ to $\p_\ii\p_\kk$,
and to the projection of $\App(\p_\ii,\p_\jj,\p_k)$, 
which is the line from $\p_\ii\p_j$ to $\p_\ii\p_\kk$.
Hence \ref{i- proj de xEi sur tree at p_i} is proven.
Assertion \ref{i- proj de xEi sur tree at D_i}
is proven in the same way.

We now prove \ref{i- proj de xEj sur tree at p_i}.
By \ref{i- proj de xEi sur tree at D_i} applied to $\vstE_\jj$, 
we have that $\projE_{\D_\jj}(\vstE_\jj)$
is  the center  of the ideal tripod 
$\p_\jj$, $\p_{\jj\kk}=\D_\jj\cap\D_\kk$, $\D_\jj\cap\D_\ii$. 
As $\vstE_\jj$ is on a geodesic from $\D_\jj$ to $\p_\ii$, 
we may deduce that $\projE_{\p_\ii}(\vstE_\jj)$
is  the center  of the ideal tripod 
$\p_\ii\p_\jj$, $\p_\ii\p_{\jj\kk}$, $\D_\ii$ 
(using the canonical isomorphim $\Es_{\D_\jj}\isomto \Es_{\p_\ii}$).
The last assertion \ref{i- proj de xEj sur tree at D_i}
has identical proof.
\end{proof}

\bibliographystyle{alpha}
\bibliography{triples}

\end{document}